\documentclass[12pt,twoside,reqno]{amsart}
\usepackage[colorlinks=true,citecolor=blue]{hyperref}
\usepackage{mathptmx, amsmath, amssymb, amsfonts, amsthm, mathptmx, enumerate, color}
\usepackage{graphicx}
\usepackage{epstopdf}

\usepackage{bm}

\newtheorem{theorem}{Theorem}[section]
\newtheorem{corollary}{Corollary}[section]

\newtheorem{proposition}{Proposition}[section]
\theoremstyle{definition}
\newtheorem{definition}{Definition}[section]
\newtheorem{example}{Example}[section]

\numberwithin{equation}{section}

\newcommand{\jbe}{\begin{equation}}
\newcommand{\jen}{\end{equation}}

\newcommand{\kg}{\mbox{}\hspace{0.3in}}

\begin{document}
\setcounter{page}{1}

\vspace*{1.0cm}
\title[Convertible Nonconvex Optimization]
{Optimization Conditions and Decomposable Algorithms for Convertible Nonconvex Optimization}
\author[M. Jiang,R. Shen, Z.Q. Meng, C.Y. Dang]{Min Jiang$^{1}$, Rui Shen$^{2,*}$, Zhiqing Meng$^{1}$, Chuangyin Dang$^{3}$}
\maketitle
\vspace*{-0.6cm}

\begin{center}
{\footnotesize {\it

$^1$School of Management, Zhejiang University of Technology,Hangzhou, Zhejiang, 310023, China\\
$^2$School of Economics, Zhejiang University of Technology,Hangzhou, Zhejiang, 310023, China\\
$^3$Department of Advanced Design and System Engineering,City University of Hong Kong, Kowloon, Hong Kong
}}\end{center}

\vskip 4mm {\small\noindent {\bf Abstract.}
This paper defines a convertible nonconvex function(CN function for short) and a weak (strong) uniform (decomposable, exact) CN function, proves the optimization conditions for their global solutions and proposes algorithms for solving the unconstrained optimization problems with the decomposable CN function. First, to illustrate the fact that some nonconvex functions, nonsmooth or discontinuous, are actually weak uniform CN functions, examples are given. The operational properties of the CN functions are proved, including addition, subtraction, multiplication, division and compound operations.  Second, optimization conditions of the global optimal solution to unconstrained optimization with a weak uniform CN function are proved.  Based on the unconstrained optimization problem with the decomposable CN function, a decomposable algorithm is proposed by its augmented Lagrangian penalty function and its convergence is proved. Numerical results show that an approximate global optimal solution to unconstrained optimization with a CN function may be obtained by the decomposable algorithms. The decomposable algorithm can effectively reduce the scale in solving the unconstrained optimization problem with the decomposable CN function. This paper provides a new idea for solving unconstrained nonconvex optimization problems.

 }

\renewcommand{\thefootnote}{}
\footnotetext{ $^*$Corresponding author.
\par
E-mail addresses: jiangmin@zjut.edu.cn (M. Jiang), shenrui126@126.com (R. Shen), mengzhiqing@zjut.edu.cn (Z.Q.Meng), mecdang@cityu.edu.hk (C.Y.Dang).
\par
Received January 16, 2022;}% Accepted February 16, 2016. }

%\subclass{90C06 \and 90C25 \and 90C26 \and 90C59}

%All acknowledgements should be placed in the back of the paper after Conclusions..

\section{Introduction}

In this paper the following unconstrained optimization (convertible nonconvex optimization, CNO)  with a weak uniform (decomposable) convertible nonconvex(CN) function is considered:
\begin{eqnarray*}
\mbox{(CNO)}\qquad & \min\; & f(\bm x) \\
& \mbox{s.t.}\; & \bm x\in R^n,
\end{eqnarray*}
where $f: R^n \rightarrow R$ is neither convex nor smooth. In machine learning, there are many nonconvex, nonsmooth, non-Lipschitz and discontinuous optimization problems in \cite{Li1,Mohri,Wang,Zhang}. So, to solve these problems, theoretical tools of nonsmooth and nonconvex functions are needed, such as the subdifferentiable, general convex, smoothing and so on in \cite{Bagirov1,Chieu,Clarke,Grant,Rockafellar}. A new nonconvex function is defined in this paper, which is called the (weak or strong uniform) CN function in Definition 2.2, where  the CN function is a nonconvex nonsmooth function form that can be transformed into a convex smooth function with convex equality constraints. The CN function somewhat relates to upper -$UC^k$ function \cite{Bougeard,Daniilidis,Hare2,Ngai,Rockafellar,Spingarn} and factorable nonconvex function \cite{Bunin,Granvilliers,He,Jackson,Lundell1,McCormick,Tawarmalani2,Serrano}.

%For example, the loss function $f$ is nonconvex nonsmooth for multi classification problem in \cite{Mohri}:
%\begin{eqnarray*} f(\bm x,\bm \theta)=\sum\limits_{i\in I}(sgn(\bm a_{i}^\top \bm x-\theta_i)-b_i)^2,\end{eqnarray*}
%where $a_i\not=0(i=1,2,\cdots,n)$ and $\bm x\in R^n$.
%where $\mbox{sgn}(\cdot)$ is a symbolic function for $\bm a_i,\bm x\in R^n$ and $(\bm a_i,b_i)(i\in I)$ is a sample data.
%A loss function for regression problem in \cite{Mohri} $f(\bm x)=\sum\limits_{i\in I}|(|\bm a_i^\top \bm x|-b_i)|$ is nonconvex nonsmooth.
%A loss function $f(\bm x)=\lambda \|\bm x\|_0+\sum\limits_{i\in I}(\bm a_i^\top \bm x-b_i)^2$ is nonconvex nonsmooth for sparse regression problem in \cite{Li1}.

The lower(upper)-$C^k$ function was suggested by Professor R. T. Rockafellar\cite{Rockafellar}. The class of lower-$C^1$ functions is first introduced by Spingarn in \cite{Spingarn}.  In his work, Spingarn
showed that these functions are (Mifflin) semi-smooth and Clarke regular and are characterized by a generalized monotonicity property of their subgradients, called submonotonicity. The definition of lower(upper)-$C^k$ function is given as follows\cite{Daniilidis}.
\begin{definition}
Let $U$ be an open subset of $R^n$ and $k\in N$.
 Function $f: U\to R^1$ is called lower-$C^k$(for short, $LC^k$), if for every $x_0\in U$ there exist $\delta> 0$, compact topological space $S$, and a jointly continuous function $F: B(x_0, \delta)\times S\to R^1$ satisfying
$$f(x) = \max\limits_{s\in S} F(x, s),\mbox{    for all  }  x\in B(x_0, \delta),$$
such that all derivatives of $F$ up to order $k$ with respect to $x$ exist and are jointly continuous. If $-f$ is lower-$C^k$, then $f$ is called upper-$C^k$ function.
\end{definition}

The lower(upper)-$C^k$ function is nonconvex or nondifferentiable, but it is locally Lipschitz approximately convex functions in \cite{Daniilidis}.
Research on the lower(upper)-$C^k$ functions is done on subdifferentiation and optimization in \cite{Hare1,Hare2,Hare3,Hare4}.  The Moreau envelopes  $erf$:
$$erf(x) := \inf\limits_{w}\{ f(w)+\frac{r}{2}|w-x|^2\} $$
is lower-$C^2$ in \cite{Bougeard,Hare2,Ngai} such that subdifferential of the lower(upper)-$C^k$ functions can solve nonconvex optimization by prox-regularity and the proximal mapping(operator) in \cite{Hare3}.  Chieu et al. proved
second-order necessary and sufficient conditions for lower-$C^2$ functions to be convex and strongly convex in \cite{Chieu}.

Some methods for non-smooth non-convex optimization programs with lower(upper)-$C^k$ functions have been studied in \cite{Dao,Hare5,Hare6,Noll}. Dao developed a nonconvex bundle method based on the downshift mechanism and a proximity control management technique to solve nonconvex nonsmooth constrained optimization problems, where he  proved its global convergence in the sense of subsequences for both classes of lower-$C^1$ and upper-$C^1$ in \cite{Dao}.
 Hare et al.  studied two  proximal bundle methods for nonsmooth nonconvex optimization in \cite{Hare5,Hare6} by proximal mapping on lower-$C^2$ functions. Noll defined a first-order model of $f$ as an extend case of lower-$C^k$ function and presented a bundle method in \cite{Noll} as follows.
 \begin{definition}
 A function $\phi: R^n\times R^n\to R^1$ is called a first-order model of $f$ on $\omega\subset
R^n$, if $\phi(\cdot, x)$ is convex for every fixed $x\in \omega$, and if the following axioms are satisfied:

(M1) $\phi(x, x) = f(x)$ and $\partial_1\phi(x, x)\subset \partial f(x)$.

(M2) For every sequence $y_j\to x$ there exists $\epsilon_j\to 0^+$ such that $f(y_j)\leq  \phi(y_j, x)+\epsilon_j\|y_j-x\|$ for all $j\in {N}$.

(M3) For sequences $y_j\to y \in R^n$ and $x_j\to x$ in $\omega$ one has $f(y_j)\leq \limsup_{j\to \infty}\leq  \phi(y, x)$ for all $j\in {N}$.
\end{definition}
Clearly, if $f$ a first-order model, $f$ is not necessarily lower-$C^k$, and the reverse is not necessarily true.

On the other hand, the branch-and-bound method in conjunction with underestimating convex problems had been proved as an effective method to solve global nonconvex optimization problems in \cite{Al-Khayyaltt,Bao1,Tawarmalani1}. Almost all the methods used to solve nonconvex optimization are to construct many convex relaxation subproblems with convex envelopes and convex underestimating, as in
\cite{Bao1,Sherali2,Serrano,Tawarmalani2}. Based on this idea, the factorable programming technique, one of the most popular approaches for constructing convex relaxations of nonconvex optimization problems including problems with convex-transformable functions, was given in \cite{McCormick}. Due to its simplicity, factorable programming technique is included in most global optimization packages such as BARON(1996), ANTIGONE(2014), etc\cite{Nohra}. But,  Nohra and Sahinidis(2018) pointed out that a main drawback of factorable programming technique is that it often results in large relaxation gaps in \cite{Nohra}.

In 1976, McCormick(1976)\cite{McCormick} first defined factorable nonconvex function, but factorable nonconvex function is not necessarily lower-$C^1$, such as $f(x)=|x|^{0.1}+|x+1|^{0.2}$ on $x\in R^1$, because $f(x)=|x|^{0.1}+|x+1|^{0.2}$ is not locally Lipschitz in \cite{Chen1}.
%Jackson and Mccormick (1988) discussed second-order sensitivity analysis for local optimal solution to factorable programming to ensure the stability of the designed algorithm\cite{Jackson}. A good merit of factorable nonconvex programs is that it can solve nonconvex discontinuous optimization problems, such as solving global optimization of mixed-integer nonlinear programs by designed branch-and-bound algorithm in \cite{Tawarmalani1}.
In fact, the factorable nonconvex functions in \cite{Jackson,Lundell1,McCormick,Tawarmalani2} may be special CN functions (see Definition 2.2). In recent years,  research on nonconvex factorable programming further shows its effectiveness in solving the global optimization, as shown in \cite{Bunin,Granvilliers,He,Serrano}.

  %When $f$ is a CN function, there are sufficient conditions that help determine the global optimal solution to (CNO) and such sufficient conditions  are easily proved. But, (CNP) with a CN function increases the scale of problem solved. Many CN functions have decomposability. (CNP) with a decomposable CN function can be decomposed into smaller optimization problems, which provides us a feasible idea of transforming the nondifferentiable nonconvex optimization problem into differentiable convex optimization.

There are many  CN functions that are not upper-$C^k$ functions or factorable  functions, such as $|x|_0$ because  upper-$C^k$ functions are continuous as per \cite{Hare1}. So, a CN function is not necessarily an upper-$C^k$ function or a factorable nonconvex function. There are three differences between factorable functions and CN functions.

 (1) They are different in the decomposition of functional representation. Each function $X^i(x)$ \emph{w.r.t} a single variable in any form of a factorable function is not necessarily a convex or concave function in \cite{McCormick}. However, each function $g_i(\bm x,\bm y)$ in any CN form of a CN function is convex  in Definition 2.2.

 (2) To estimate the factorable function, it is necessary to underestimate/ overestimate the  convex/concave functions of each $X^i(x)$, while CN function does not require estimation of its convex envelope(underestimating) function.

(3) The method of solving optimization problems with a CN function differs from the method of solving nonconvex factorable programming. Factorial programming  solves its relaxation problem to obtain an approximate global optimal solution. To solve the approximate global optimal solution to optimization problem with a CN function, its equivalent optimization is needed. %Our idea is  as follows.

In order to solve (CNO), Jiang et al(2021) \cite{Jiang} have discussed optimal conditions, Lagrangian dual and an algorithm for the unconstrained CN  optimization problems. Different from that, in this paper, a weak uniform  CN function is defined, and weak uniform and decomposable  weak uniform of CN function, optimization conditions and decomposable algorithms for (CNO) with weak uniform  CN function are studied. The main contributions of this paper are as follows: (1) a weak uniform  CN function is proposed, (2) the sufficient conditions of the optimal solution to optimization problems with such CN function are proved, and (3) a decomposable algorithm for the optimization problem with the decomposable CN function is proposed. This paper provides a new method to solve the difficulties in solving nonconvex or nonsmooth optimization problems.

The remainder of the paper is organized as follows. In Section 2,  a CN function and a weak (strongy) uniform (decomposable, exact) CN function are defined respectively. Some examples are given. The operational properties of the CN functions are given.  In section 3, optimization conditions of the global optimal solution to unconstrained optimization with a weak uniform CN function are proved. In Section 4, for the decomposable CN function, a decomposable algorithm is proposed by its augmented Lagrangian penalty function and its convergence is proved. In section 5, the conclusion is given.

\section{Weak Uniform CN Function}

In this section, a (weak, strong uniform, exact) CN function is defined. Some examples are given to show that some nonconvex or discontinuous functions are differentiable (weak uniform) CN ones.

\begin{definition}%{\bf Definition 1}
Let function $g:R^n\times R^m\to R^1$ be differentiable. For all $\bm d\in R^n\times R^m$ and all $(\bm x,\bm y)\in R^n\times R^m$, if there is a positive semidefinite matrix $B(\bm x,\bm y)$ such that
\begin{eqnarray}
    g((\bm x,\bm y)+\bm d)-g(\bm x,\bm y)\geq \nabla g(\bm x,\bm y)^\top \bm d+ \frac{1}{2} \bm d^\top B(\bm x,\bm y) \bm d,\label{d1}
\end{eqnarray}
then $g$ is called a weak uniform convex function which has two cases as follows.

(i) If there is a positive definite matrix $B(\bm x,\bm y)$ such that (\ref{d1}) holds
then $g$ is called a strong uniform convex function.

(ii) If there is an $\bar{\rho}>0$ such that
\begin{eqnarray}
    g((\bm x,\bm y)+\bm d)-g(\bm x,\bm y)\geq \nabla g(\bm x,\bm y)^\top \bm d+\frac{\bar{\rho}}{2}\|\bm d\|^2,\label{d2}
\end{eqnarray}
then $g$ is called a uniform convex function.
\end{definition}

It is clear that a uniform convex function $g$ is not only a weak uniform convex function but also a strong uniform convex function. A strong uniform convex function is a weak uniform convex function.
And a weak uniform convex function is a convex function, as shown in the following examples.

%{\bf Example 1}
\begin{example} $g(x,y)=(x+y-1)^2$ for $(x,y)\in R^1\times R^1$
 is a weak uniform convex function. But, $g(x,y)$ is not a strong uniform convex function.
\end{example}
\begin{example}
%{\bf Example 2}
$g(x,y)=x^4+y^4$ for $(x,y)\in R^1\times R^1$ is a strong uniform convex function for $(x,y)\not=0$. But, $g(x,y)$ is not a uniform convex function.
\end{example}
\begin{example}
%{\bf Example 3}
$g(\bm x,\bm y)=(\bm x,\bm y)^\top A (\bm x,\bm y)+\bm c^\top(\bm x,\bm y)$ is a weak uniform convex function, where $A$ is a positive semidefinite matrix, $(\bm x,\bm y)\in R^n\times R^m$.
\end{example}

We have the following conclusion.

\begin{proposition}%{\bf Proposition 1}
Let function $g:R^n\times R^m\to R^1$ be twice continuously differentiable and matrix $B(\bm x,\bm y)$  be given for $(\bm x,\bm y)\in R^n\times R^m$.  Then, $g(\bm x,\bm y)$ is a weak(strong) uniform convex function if and only if
\begin{eqnarray*}
 \bm d^\top \nabla^2 g(\bm x,\bm y)\bm d\geq  \bm d^\top B(\bm x,\bm y) \bm d\geq(>) 0,\ \forall (\bm x,\bm y),\forall \bm d\in R^n\times R^m.
\end{eqnarray*}
% (ii) $g(\bm x,\bm y)$ is a strong uniform convex function if and only if
%\begin{eqnarray*} \bm d^\top \nabla^2 g(\bm x,\bm y)\bm d\geq \bm d^\top B(\bm x,\bm y) \bm d> 0,\ \forall (\bm x,\bm y),\forall \bm d\in R^n\times R^m\backslash\{\bm 0\};\end{eqnarray*}
%(iii) $g(\bm x,\bm y)$ is a  uniform convex function if and only if
%\begin{eqnarray*} \bm d^\top \nabla^2 g(\bm x,\bm y)\bm d\geq \bar{\rho}\|\bm d\|^2> 0,\ \forall (\bm x,\bm y), \forall \bm d\in R^n\times R^m\backslash\{\bm 0\}.\end{eqnarray*}
\end{proposition}

Now based on the definition of  weak uniform convex function, let's define the weak uniform CN function, strong uniform CN function and uniform CN function.

%{\bf Definition 2}
\begin{definition} Let $S=S_1\times S_2\subset R^n\times R^m$ be a convex set.
 Let functions $g_i:R^n\times R^m\to R^1$ ($i=1,2,\cdots,r$ and $r\geq 1$) be convex on $S$ and $g: R^n\times R^m\to R^1$ be a convex function on $S$.  Let function $f: R^n \rightarrow R^1$ be nonconvex. Let
$$\bm g(\bm x,\bm y)=(g_1(\bm x,\bm y),g_2(\bm x,\bm y),\cdots,g_r(\bm x,\bm y))^\top.$$
Let a set
\begin{eqnarray}
  X(\bm g)=\{(\bm x,\bm y)\in S\mid g_i(\bm x,\bm y)=0,\  i=1,2,\cdots,r\}. \label{d3}
\end{eqnarray}
If for each $\bm x\in S$ and  $f(\bm x)$ there is a $\bm y\in R^m$ such that $(\bm x,\bm y)\in X(\bm g)$ and
\begin{eqnarray*}
f(\bm x)=g(\bm x,\bm y)=\min\limits_{(\bm x,\bm y')\in X(\bm g)} g(\bm x,\bm y'),\label{d4}
\end{eqnarray*}
then $f$ is called a convertible nonconvex(CN) function on $S$ (when $S=R^n\times R^m$, the term "on $S$" is omitted).  $[g:g_1,g_2,\cdots,g_r]$ is called a convertible nonconvex(CN) form of $f$ on $S$, briefing as $f=[g:g_1,g_2,\cdots,g_r]$. For $f$, the number of its  CN form is more than one. In particular,

(i) if $g$ is a weak uniform convex function on $S$,  then $f$ is called a weak uniform CN function on $S$;

(ii) if $g$ is a strong uniform convex function on $S$, then  $f$ is called a strong uniform CN function on $S$;

(iii) if $g$ is a uniform convex function on $S$, then $f$ is called a uniform CN function on $S$.

Furthermore, if for each $\bm x\in S$ and  $f(\bm x)$ there is a $\bm y\in R^m$ such that $(\bm x,\bm y)\in X(\bm g)$ and
\begin{eqnarray*}
f(\bm x)=g(\bm x,\bm y)= g(\bm x,\bm y'), \ \ \ \forall(\bm x,\bm y')\in X(\bm g),\label{d4}
\end{eqnarray*}
then $f$ is called an exact convertible nonconvex(CN) function on $S$.  $[g:g_1,g_2,\cdots,g_r]$ is called an exact convertible nonconvex(CN) form of $f$ on $S$. Then, $f$ is called a (weak or strong) uniform exact CN function on $S$ in (i),(ii) or (iii) if the condition of in (i),(ii) or (iii) holds accordingly.
\end{definition}

When $S=R^n\times R^m$, the term "on $S$" above is all omitted. It is clear that a uniform CN function is a CN function, a weak uniform CN function and a strong uniform CN function. A weak uniform convertible convex function is a CN function.

{\bf Remark } It is clear that the definition of  CN function differs form that of upper-$C^k$ function  \cite{Daniilidis} and a first-order model \cite{Noll}. For example, 0-norm $\|\bm x\|_0$ is not an upper-$C^k$ function  \cite{Daniilidis,Hare1} or a first-order model, because it is not continuous. But, $\|\bm x\|_0$ is a CN function in Example 2.7. Hence, CN function contains a wider range of functions than upper-$C^k$ functions  and first-order model.

Definition 2.2 means that a nonconvex function $f(\bm x)$ may be converted into a (weak uniform) convertible convex function. By Definition 2.2, a set is defined by
\begin{eqnarray}
X(f)=\{(\bm x,\bm y)| f(\bm x)=g(\bm x,\bm y), \forall (\bm x,\bm y)\in X(\bm g)\}.\label{d4}
\end{eqnarray}

For a fixed $(\bm x,\bm y)$, two sets are defined by
\begin{eqnarray*}
Y_g(\bm x)=\{\bm y\in R^m| (\bm x,\bm y)\in X(\bm g)\},
X_g(\bm y)=\{\bm x\in R^m| (\bm x,\bm y)\in X(\bm g)\}.
\end{eqnarray*}

The following conclusion is clear.

%{\bf Proposition 2}
\begin{proposition}
Let $f$ be a CN function on $S$. Then, (1) $f$ is an exact CN function on $S$ if and only if $X(f)=X(\bm g)$; (2) $X(f)$ and $X(\bm g)$ are closed sets for $S=R^n$; (3) $f(\bm x)=\min\limits_{\bm y\in Y_g(\bm x)} g(\bm x,\bm y)$; and (4)  $\min\limits_{\bm x\in X_g(\bm y)} f(\bm x)\leq \min\limits_{\bm x\in X_g(\bm y)} g(\bm x,\bm y)$.
\end{proposition}

Proposition 2.2 shows that if $f$ is a CN function on $S$, $f(\bm x)=g(\bm x,\bm y)\leq g(\bm x,\bm y')$ for all $(\bm x,\bm y)\in X(f)$ and $(\bm x,\bm y')\in X(\bm g)\backslash X(f)$. $(\bm x,\bm y)\in X(f)$ is called a CN point of $f$.

Next,  an  example is given to show that the number of weak uniform (exact) CN forms could be more than one.

\begin{example}
%{\bf Example 4}
Non-convex function $f(x_1,x_2)=2x_1x_2$ is a weak uniform exact CN function. One of its weak uniform exact CN forms is $
f=[(x_1+x_2)^2-y_1-y_2:x_1^2-y_1,x_2^2-y_2].$
A second weak uniform exact CN form of $f(x_1,x_2)=2x_1x_2$ is $f=[0.5(x_1+x_2)^2-0.5y_1:x_1^2+x_2^2-y_1]$.
Hence, it is understood that there are more than one weak uniform exact CN form.
\end{example}
%\begin{example}

Next,  some operational properties of the weak uniform (exact) CN function are easily proved as follows.

\begin{proposition}%{\bf Proposition 3}
If $f:R^n\to R$ is an exact CN function on $S$, then
$-f$ is an exact CN function  on $S$.
\end{proposition}

\begin{proposition}%{\bf Proposition 4}
 If $f_1,f_2:R^n\to R$ are (weak uniform, strong uniform)(exact) CN functions  on $S$, then
$\alpha_1 f_1+\alpha_2 f_2$ is a (weak uniform, strong uniform) (exact) CN function  on $S$ for any $\alpha_1,\alpha_2>0$.
Especially, if $f_1,f_2:R^n\to R$ are exact CN functions  on $S$, then
$\alpha_1 f_1+\alpha_2 f_2$ is an exact CN function  on $S$ for any $\alpha_1,\alpha_2\in R^1$.

\end{proposition}

\begin{proposition}%{\bf Proposition 5}
 If $f_1,f_2:R^n\to R$ are (weak uniform or strong uniform) exact CN functions  on $S$, then
$f_1f_2$ is a weak uniform exact CN function  on $S$.

%(ii) If $f_1,f_2:R^n\to R$ are weak uniform exact CN functions  on $S$, then $f_1f_2$ is a weak uniform exact CN function  on $S$.

%(iii) If $f_1,f_2:R^n\to R$ are strong uniform exact CN functions, then $f_1f_2$ is a weak uniform exact CN function  on $S$.
\end{proposition}

\begin{proposition}%{\bf Proposition 6}
 If $f_1,f_2:R^n\to R$ are (weak uniform, strong uniform) exact CN functions  on $S$, then
$\frac{f_1}{f_2}$ is a weak uniform exact CN function  on $S$.

\end{proposition}

%Then, by the proof of (i), (ii) and (iii) hold.\\

\begin{proposition}%{\bf Proposition 7}
 If $f:R^n\to R$ is a (weak uniform, strong uniform)(exact) CN function  on $S$ and $\phi:R\to R$ is a monotone increasing convex function, then $\phi(f)$ is a (weak uniform, strong uniform)(exact) CN function  on $S$.

\end{proposition}

\begin{example}%{\bf Example 5}
Let a weak uniform DC function $f(\bm x,\bm y)=d(\bm x,\bm y)-c(\bm x,\bm y)$ for $(\bm x,\bm y)\in R^n\times R^m$, where  $d(\bm x,\bm y)$ and $c(\bm x,\bm y)$ are weak uniform convex functions on $(\bm x,\bm y)\in R^n\times R^m$. Let $g(\bm x,\bm y,z)=d(\bm x,\bm y)-z$ and $g_1(\bm x,\bm y,z)=c(\bm x,\bm y)-z$.
So, $f(x,y)$ is a weak uniform CN function.
\end{example}

 By Proposition 2.3-2.7, some polynomial functions are CN. For example, multi-convex function $f(\bm x)=x_1x_2\cdots x_n$ is a weak uniform CN by Proposition 2.4. Therefore, weak uniform CN functions cover a wide range of non-convex functions. To illustrate, some examples are given as follows.

\begin{example}%{\bf Example 6}
(Example 2.1 in  \cite{Chen1}) A CN form of nonsmooth function $f(x_1,x_2)=(x_1+x_2-1)^2+\lambda(|x_1|^\frac{1}{2}+|x_2|^\frac{1}{2})$ is defined as
 \begin{eqnarray*}
[(x_1+x_2-1)^2+\lambda(y_1+y_4):y_1^4-y_3,x_1^2-y_3,y_2^2-y_1,y_4^4-y_6,x_2^2-y_6,y_5^2-y_4],
 \end{eqnarray*}
where $\lambda>0$. So, $f(x)$ is a weak uniform exact CN function.
\end{example}
\begin{example}%{\bf Example 7.}
 A CN form of the nonconvex and discontinuous function $f(x_1,x_2)$ $=(x_1+x_2-1)^2+\lambda\|(x_1,x_2)\|_0$ is defined as
 \begin{eqnarray*}
 f&=&[(x_1+x_2-1)^2+\lambda(y_1^2+y_3^2):(x_1+y_1-1)^2-y_2,x_1^2+(y_1-1)^2-y_2,\\&& y_1^2-y_1,
 (x_2+y_3-1)^2-y_4,x_2^2+(y_3-1)^2-y_4,y_3^2-y_3],
\end{eqnarray*}
 where $\lambda>0$. So, $f(x)$ is a weak uniform CN function, but is not exact.
\end{example}

The above examples show that some nonsmooth, nonconvex or discontinuous functions are CN functions as shown by their twice differentiable and weak uniform CN forms.
%, for example $ h(x)+\lambda\|x\|_0$ in  a sparse optimization problem, where $\|x\|_0$ is 0-norm number and $h(x)$ is convex. It is easy to know that there are many of nonsmooth, nonconvex or discontinuous functions are CN functions by  Proposition 2.3-2.7, because their CN forms may take many simple forms of CN functions on $S=R^1$, such as $\|x\|_0$,$\|x\|^a$, $x^a$, etc.where $0<a<1$. Hence, anyone polynomial function is a CN function and its exact CN form is easily obtained
  By Proposition 2.3-2.7, in the fields of machine learning, image processing and signal processing, many nonconvex functions are combinations of simple CN functions, for example, in  \cite{Nohra} the examples of nonconvex functions:
\begin{eqnarray*}
\prod\limits_{i=1}^m x_i^{\alpha_i},\
f(\bm x)^\alpha g(\bm x)^\beta, \
\frac{f(\bm x)^\alpha exp(g_0(\bm x))}{k_0+k_1\sum\limits_{i=1}^m exp(g_i(\bm x))},
\end{eqnarray*}
where $f(\bm x),g(\bm x),g_i(\bm x)(i=0,1,\cdots,m)$ are convex on $\bm x$, are CN functions.
The compiled functions show themselves in a wide variety of scientific and engineering applications. Their exact CN forms are easily obtained  by Proposition 2.3-2.7.

%DC function is a very important class of nonconvex functions in THi et al(2018). All DC functions are CN. Proposition 8 tells us that all second-order continuously differentiable functions on $R^n$ are CN.

%{\bf Proposition 8} If $f:R^n\to R^1$ is a second-order continuously differentiable function on $R^n$, then $f(x)$ is CN.

%{\bf Proof.} By Hartman (1959), if $f:R^n\to R^1$ is a second-order continuously differentiable function on $R^n$, $f(x)$ is a DC function. Hence, $f(x)$ is a CN function.\\

\section{Optimization Conditions of (CNO)}

In this section, it is assumed that $f$ is a weak uniform  CN function or a CN function with $f=[g:g_1,g_2,\cdots,g_r]$, where  $g,g_1,g_2,\cdots,g_r$ are twice differentiable.
But, $f(\bm x)$ is not necessarily differentiable on $\bm x\in R^n$.

With the (weak uniform) CN form of $f$,
consider the following constrained optimization problem:
\begin{eqnarray*}
\mbox{(CNP)}\qquad & \min\limits_{\bm x,\bm y)\in R^n\times R^m} & g(\bm x,\bm y) \\
& \mbox{s.t.}\; & g_i(\bm x,\bm y)=0,i=1,2,\cdots,r.
\end{eqnarray*}
It is clear
$\min\limits_{\bm x\in R^n} \ f(\bm x)=\min\limits_{(\bm x,\bm y)\in X(f)} \ g(\bm x,\bm y)=\min\limits_{(\bm x,\bm y)\in X(\bm g)} \ g(\bm x,\bm y).$

Let a directional sets at a fixed $(\bm x,\bm y)$ be defined by
\begin{eqnarray*}
T(\bm x,\bm y)=\{\bm d\in R^n\times R^m&|& \nabla  g_i(\bm x,\bm y)^\top\bm d\leq 0,i=1,2,\cdots,r\} \label{d9}
\end{eqnarray*}

Let us prove the sufficient condition of an optimal solution to (CNO) by solving (CNP).

\begin{theorem}%{\bf Theorem 1}
Suppose that $(\bm x^*,\bm y^*)\in X(f)$ and  $f$ is a weak uniform CN function. Let the problem
\begin{eqnarray*}
\mbox{(WCNP)}(\bm x^*,\bm y^*)\qquad & \min\; & \nabla g(\bm x^*,\bm y^*)^\top \bm d +\frac{1}{2} \bm d^\top B(\bm x^*,\bm y^*) \bm d  \\
& \mbox{s.t.}\; &  \bm d\in T(\bm x^*,\bm y^*).
\end{eqnarray*}
  If $\bm d^*$  is an optimal solution  to (WCNP)$(\bm x^*,\bm y^*)$ such that $\nabla g(\bm x^*,\bm y^*)^\top \bm d^*+\frac{1}{2}  \bm d^{*T} B(\bm x^*,\bm y^*) \bm d^* \geq 0$, then $(\bm x^*,\bm y^*)$  is an optimal solution  to (CNP) and $\bm x^*$ is an optimal solution to (CNO). Furthermore, if there is a $\bm d'$ such that $\nabla g_i(\bm x^*,\bm y^*)^\top  \bm d'< 0,i=1,2,\cdots,r,$ then there are $\alpha_i^*\geq 0, i=1,2,\cdots,r$ such that
  \begin{eqnarray}
 B(\bm x^*,\bm y^*)\bm d^* +\nabla g(\bm x^*,\bm y^*) + \sum\limits_{i=1}^r \alpha_i^* \nabla g_i(\bm x^*,\bm y^*)=0, \label{d21}\\
  \alpha_i^*\nabla g_i(\bm x^*,\bm y^*)^\top  \bm d^*= 0,i=1,2,\cdots,r. \label{d22}
\end{eqnarray}
Conversely, if there are $\bm d^*\in T(\bm x^*,\bm y^*)$ and $\alpha_i^*\geq 0, i=1,2,\cdots,r$ such that
(\eqref{d21}) and (\eqref{d22}) hold, then $\bm d^*$  is an optimal solution  to (WCNP)$(\bm x^*,\bm y^*)$.
\end{theorem}
  {\it Proof.} For any $(\bm x,\bm y)\in X(\bm g)$, let $\bm d=[(\bm x,\bm y)-(\bm x^*,\bm y^*)]$. Then we have
  \begin{eqnarray*}
 0=g_i(\bm x,\bm y)-g_i(\bm x^*,\bm y^*)&\geq& \nabla g_i(\bm x^*,\bm y^*)^\top[(\bm x,\bm y)-(\bm x^*,\bm y^*)],\kg  i=1,2,\cdots,r.
\end{eqnarray*}
So,  $(\bm x,\bm y)-(\bm x^*,\bm y^*)$ is a feasible solution to (WCNP)$(\bm x^*,\bm y^*)$. Then
\begin{eqnarray*}
g(\bm x,\bm y)-g(\bm x^*,\bm y^*)&\geq& \nabla g(\bm x^*,\bm y^*)^\top\bm d+\frac{1}{2} \bm d^\top B(\bm x^*,\bm y^*) \bm d ,\\
&\geq &\nabla g(\bm x^*,\bm y^*)^\top \bm d^*+\frac{1}{2} \bm d^{*\top} B(\bm x^*,\bm y^*) \bm d^* \geq 0.
\end{eqnarray*}
Hence, $(\bm x^*,\bm y^*)$ is an optimal solution to (CNP) and $\bm x^*$ is an optimal solution to (CNO).
Since (WCNP)$(\bm x^*,\bm y^*)$  is a convex programm, by KKT condition, the second conclusion (\eqref{d21}) and (\eqref{d22}) are true. Conversely, the conclusion is clear.

 By  Theorem 3.1, we have three corollaries as follows.

\begin{corollary} %{\bf Corollary 1}
 Suppose that $(\bm x^*,\bm y^*)\in X(f)$ and  $f$ is a uniform CN function. Let
\begin{eqnarray*}
\mbox{(UCNP)}(\bm x^*,\bm y^*)\qquad & \min\; & \nabla g(\bm x^*,\bm y^*)^\top \bm d +\frac{\bar{\rho}}{2}\|\bm d\|^2 \\
& \mbox{s.t.}\; &  \bm d\in T(\bm x^*,\bm y^*).
\end{eqnarray*}
 If $\bm d^*$  is an optimal solution  to (UCNP)$(\bm x^*,\bm y^*)$ such that $\nabla g(\bm x^*,\bm y^*)^\top \bm d^*+\frac{\bar{\rho}}{2}\|\bm d^*\|^2 \geq 0$, then $(\bm x^*,\bm y^*)$  is an optimal solution  to (CNP) and $\bm x^*$ is an optimal solution to (CNO). Furthermore, if there is a $\bm d'$ such that $\nabla g_i(\bm x^*,\bm y^*)^\top  \bm d'< 0,i=1,2,\cdots,r,$ then there are $\alpha_i^*\geq 0, i=1,2,\cdots,r$ such that
  \begin{eqnarray}
 \bar{\rho}\bm d^* +\nabla g(\bm x^*,\bm y^*) + \sum\limits_{i=1}^r \alpha_i^* \nabla g_i(\bm x^*,\bm y^*)=0, \label{d23}\\
  \alpha_i^*\nabla g_i(\bm x^*,\bm y^*)^\top  \bm d^*= 0,i=1,2,\cdots,r. \label{d24}
\end{eqnarray}
Conversely, if there are $\bm d^*\in T(\bm x^*,\bm y^*)$ and $\alpha_i^*\geq 0, i=1,2,\cdots,r$ such that
(\eqref{d23}) and (\eqref{d24}) hold, then $\bm d^*$  is an optimal solution  to (UCNP)$(\bm x^*,\bm y^*)$.
\end{corollary}

\begin{corollary}% {\bf Corollary 2}
Suppose that $(\bm x^*,\bm y^*)\in X(f)$ and $f$ is a CN function.
\begin{eqnarray*}
\mbox{(LCNP)}(\bm x^*,\bm y^*)\qquad & \min\; & \nabla g(\bm x^*,\bm y^*)^\top \bm d  \\
& \mbox{s.t.}\; & \bm d\in T(\bm x^*,\bm y^*).
\end{eqnarray*}
 If $\bm d^*$  is an optimal solution  to (LCNP)$(\bm x^*,\bm y^*)$ such that $\nabla g(\bm x^*,\bm y^*)^\top \bm d^*\geq 0$, then $(\bm x^*,\bm y^*)$  is an optimal solution  to (CNP) and $\bm x^*$ is an optimal solution to (CNO). Furthermore, if there is a $\bm d'$ such that $\nabla g_i(\bm x^*,\bm y^*)^\top  \bm d'< 0,i=1,2,\cdots,r,$ then there are $\alpha_i^*\geq 0, i=1,2,\cdots,r$ such that
  \begin{eqnarray}
 \nabla g(\bm x^*,\bm y^*) + \sum\limits_{i=1}^r \alpha_i^* \nabla g_i(\bm x^*,\bm y^*)=0. \label{d25}\\
  \alpha_i^*\nabla g_i(\bm x^*,\bm y^*)^\top  \bm d^*= 0,i=1,2,\cdots,r. \label{d26}
\end{eqnarray}
Conversely, if there are $\bm d^*\in T(\bm x^*,\bm y^*)$ and $\alpha_i^*\geq 0, i=1,2,\cdots,r$ such that
(\eqref{d25}) and (\eqref{d26}) hold, then $(\bm x^*,\bm y^*)$  is an optimal solution  to (CNP) and $\bm x^*$ is an optimal solution to (CNO).
\end{corollary}

\begin{corollary}% {\bf Corollary 3}
Suppose that $(\bm x^*,\bm y^*)\in X(f)$ and $f$  is a CN function.
If $\nabla g(\bm x^*,\bm y^*)=0$, then $(\bm x^*,\bm y^*)$  is an optimal solution  to (CNP) and $\bm x^*$ is an optimal solution to (CNO).\\
\end{corollary}

Define a set
\begin{eqnarray*}
T_1(\bm x^*,\bm y^*)=\{ (\bm x,\bm y) \mid \nabla g_i(\bm x^*,\bm y^*)^\top  [(\bm x,\bm y)-(\bm x^*,\bm y^*)]\leq 0,i=1,2,\cdots,r.\}
\end{eqnarray*}
We have $X(\bm g)\subset T_1(\bm x^*,\bm y^*)$ and $T_1(\bm x^*,\bm y^*)-(\bm x^*,\bm y^*)=T(\bm x^*,\bm y^*)$.

By \cite{Frank}, if the objective value $G_1(\bm d)$ of (WCNP)$(\bm x^*,\bm y^*)$ has a lower bound, then there is a global optimal solution  to (WCNP)$(\bm x^*,\bm y^*)$. Then we have the following conclusion.

\begin{theorem}%{\bf Theorem 2}
Suppose that $(\bm x^*,\bm y^*)\in X(f)$ and $f$ is a weak uniform CN function.  Let
\begin{eqnarray*}
\mbox{(WCNP1)}(\bm x^*,\bm y^*)\qquad & \min\; & G_1(\bm x,\bm y)=\nabla g(\bm x^*,\bm y^*)^\top \bm [(\bm x,\bm y)-(\bm x^*,\bm y^*)] \\ && +\frac{1}{2} \bm [(\bm x,\bm y)-(\bm x^*,\bm y^*)]^\top B(\bm x^*,\bm y^*) \bm [(\bm x,\bm y)-(\bm x^*,\bm y^*)] \\
& \mbox{s.t.}\; &  (\bm x,\bm y)\in T_1(\bm x^*,\bm y^*).
\end{eqnarray*}
If the objective value $G_1(\bm x,\bm y)$ of (WCNP1)$(\bm x^*,\bm y^*)$ has a lower bound, then there is an optimal solution $( \bar{\bm x}, \bar{\bm y})$ to (WCNP1)$(\bm x^*,\bm y^*)$. Furthermore, if $G_1( \bar{\bm x}, \bar{\bm y})$ $\geq 0$, then $(\bm x^*,\bm y^*)$  is an optimal solution  to (CNP) and $\bm x^*$ is an optimal solution to (CNO), or if $G_1(\bar{\bm x},\bar{\bm y})< 0$, then $G_1(\bar{\bm x}, \bar{\bm y})=-\frac{1}{2}\bar{\bm d}^\top B(\bm x^*,\bm y^*) \bar{\bm d}$, where $\bar{\bm d}=( \bar{\bm x}, \bar{\bm y})-(\bm x^*,\bm y^*)$. If  there is $(\bm x',\bm y')$ such that $ \nabla g_i(\bm x^*,\bm y^*)^\top  [(\bm x,\bm y)-(\bm x^*,\bm y^*)]< 0,i=1,2,\cdots,r$ holds, then there are $\alpha_i^*\geq 0, i=1,2,\cdots,r$ such that
  \begin{eqnarray*}
  \bar{\bm d}^\top B(\bm x^*,\bm y^*) +\nabla g(\bm x^*,\bm y^*) + \sum\limits_{i=1}^r \alpha_i^* \nabla g_i(\bm x^*,\bm y^*)=0.
\end{eqnarray*}
\end{theorem}

%{\it Proof.} By \cite{Frank}, if the objective value $G_1(\bm x,\bm y)$ of (WCNP1)$(\bm x^*,\bm y^*)$ has a lower bound, there is a global optimal solution $( \bar{\bm x}, \bar{\bm y})$ to (WCNP1)$(\bm x^*,\bm y^*)$. So, let $\bar{\bm d}=( \bar{\bm x}, \bar{\bm y})-(\bm x^*,\bm y^*)$. Suppose that $G_1(\bar{\bm d})< 0$, then $\bar{\bm d}\not=0$,$ \nabla g(\bm x^*,\bm y^*)^\top \bar{\bm d}<0$ and $\frac{1}{2}  \bar{\bm d}^\top B(\bm x^*,\bm y^*) \bar{\bm d}>0$. Here, if $t\in (0,+\infty)$, we have $t\bar{\bm d}\in T(\bm x^*,\bm y^*)$. Hence, $G_1(\bar{\bm d})\leq G_1(t\bar{\bm d})$ holds for $t\in (0,+\infty)$. On the other hand, the problem
%$$ \min\limits_{t>0} G_1(t\bar{\bm d})=(\nabla g(\bm x^*,\bm y^*)^\top \bar{\bm d})t+(\frac{1}{2}   \bar{\bm d}^\top B(\bm x^*,\bm y^*)\bar{\bm d})t^2 $$
%has an optimal solution $t^*=-\frac{\nabla g(\bm x^*,\bm y^*)^\top \bar{\bm d}}{\bar{\bm d}^\top B(\bm x^*,\bm y^*)  \bar{\bm d}}>0$ such that
%  \begin{eqnarray*} G_1(t^*\bar{\bm d})\leq G_1(\bar{\bm d}).  \end{eqnarray*}
%So, $G_1(t^*\bar{\bm d})= G_1(\bar{\bm d})$. Since $G_1(t\bar{\bm d})$ is strictly monotonous on $t\in (0,+\infty)$, we have $t^*=1$. Then, ${\nabla g(\bm x^*,\bm y^*)^\top \bar{\bm d}}+{ \bar{\bm d}^\top B(\bm x^*,\bm y^*) \bar{\bm d}}=0$. \\

These conclusions of Theorem 3.1 and  Theorem 3.2 are true when $f$ is a strong uniform CN function.

%Theorem 3.2  shows that $(\bm x^*,\bm y^*)$  is an optimal solution  to (CNP) if $B(\bm x^*,\bm y^*)$ is positive definite or $\bm d^{*T} B(\bm x^*,\bm y^*) \bm d^*> 0$. If $(\bm x^*,\bm y^*)$  is not an optimal solution  to (CNP), $B(\bm x^*,\bm y^*)$ is not positive definite or $\bm d^{*T} B(\bm x^*,\bm y^*) \bm d^*=0$. If there is an optimal solution to (CNO), there is an optimal solution to (CNP). If there is no optimal solution to (CNO), there is no optimal solution to (CNP).

Theorem 3.2  shows that $(\bm x^*,\bm y^*)$  is an optimal solution  to (CNP) if $B(\bm x^*,\bm y^*)$ is positive definite or $\bm d^{*T} B(\bm x^*,\bm y^*) \bm d^*> 0$. If $(\bm x^*,\bm y^*)$  is not an optimal solution  to (CNP), $B(\bm x^*,\bm y^*)$ is not positive definite or $\bm d^{*T} B(\bm x^*,\bm y^*) \bm d^*=0$.
%If there is an optimal solution to (CNO), there is an optimal solution to (CNP). If there is no optimal solution to (CNO), there is no optimal solution to (CNP).

Next, another sufficient condition of the optimal solution to (CNO) is obtained by solving (CNP).

\begin{theorem}%{\bf Theorem 3}
Suppose that $(\bm x^*,\bm y^*)\in X(f)$ and  $f$ is a weak uniform CN function. Define a set
\begin{eqnarray*}
 K_w(\bm x^*,\bm y^*)&=&\{(\bm x,\bm y) \in R^n\times R^m \mid  \nabla g(\bm x^*,\bm y^*)^\top [(\bm x,\bm y)-(\bm x^*,\bm y^*)] \nonumber\\
 && +\frac{1}{2} [(\bm x,\bm y)-(\bm x^*,\bm y^*)]^\top B(\bm x^*,\bm y^*) [(\bm x,\bm y)-(\bm x^*,\bm y^*)]<0 \}. \label{d11}
\end{eqnarray*} If
%\begin{eqnarray}
$X(\bm g)\cap K_w(\bm x^*,\bm y^*)=\emptyset$ %\label{d27}
%\end{eqnarray}
holds, then $(\bm x^*,\bm y^*)$  is an optimal solution  to (CNP) and $\bm x^*$ is an optimal solution to (CNO).
\end{theorem}

%{\bf Proof.} If there is an $(\bm x',\bm y')\in X(\bm g)$ such that $g(\bm x^*,\bm y^*)> g(\bm x',\bm y')$. We have \begin{eqnarray*} 0>g(\bm x',\bm y')-g(\bm x^*,\bm y^*)&\geq& \nabla g(\bm x^*,\bm y^*)^\top\bm [(\bm x',\bm y')-(\bm x^*,\bm y^*)]\\ && +\frac{1}{2} [(\bm x',\bm y')-(\bm x^*,\bm y^*)]^\top B(\bm x^*,\bm y^*) [(\bm x',\bm y')-(\bm x^*,\bm y^*)].\end{eqnarray*} So, $(\bm x',\bm y')\in K_w(\bm x^*,\bm y^*)$. Hence, the conclusion holds.\\

By Theorem 3.3, if $\bm (\bm x^*,\bm y^*)$  is an optimal solution  to
\begin{eqnarray*}
\mbox{(WCNPP)}(\bm x^*,\bm y^*)\qquad & \min\; & \nabla g(\bm x^*,\bm y^*)^\top \bm [(\bm x,\bm y)-(\bm x^*,\bm y^*)] \\ && +\frac{1}{2} \bm [(\bm x,\bm y)-(\bm x^*,\bm y^*)]^\top B(\bm x^*,\bm y^*) \bm [(\bm x,\bm y)-(\bm x^*,\bm y^*)] \\
& \mbox{s.t.}\; &  (\bm x,\bm y)\in T_1(\bm x^*,\bm y^*)\cap X(\bm g),
\end{eqnarray*}
 then $(\bm x^*,\bm y^*)$  is an optimal solution  to (CNP) and $\bm x^*$ is an optimal solution to (CNO).

\begin{corollary}%{\bf Corollary 4}
Suppose that $(\bm x^*,\bm y^*)\in X(f)$ and  $f$ is a uniform CN function. Define a set
\begin{eqnarray}
 K_u(\bm x^*,\bm y^*)=\{(\bm x,\bm y) \in R^n\times R^m \mid \nabla g(\bm x^*,\bm y^*)^\top [(\bm x,\bm y)-(\bm x^*,\bm y^*)]\nonumber\\
 +\frac{\bar{\rho}}{2} \|(\bm x-\bm x^*,\bm y-\bm y^*)\|^2<0 \}. \label{d28}
\end{eqnarray}
 If
%\begin{eqnarray}
$X(\bm g)\cap K_u(\bm x^*,\bm y^*)=\emptyset$ % \label{d29}
%\end{eqnarray}
holds, then $(\bm x^*,\bm y^*)$  is an optimal solution  to (CNP) and $\bm x^*$ is an optimal solution to (CNO).
\end{corollary}

By Corollary 3.1 and Corollary 3.4, if $\bm (\bm x^*,\bm y^*)$  is an optimal solution  to
\begin{eqnarray*}
\mbox{(UCNP)}(\bm x^*,\bm y^*)\qquad & \min\; & \nabla g(\bm x^*,\bm y^*)^\top \bm [(\bm x,\bm y)-(\bm x^*,\bm y^*)] +\frac{\bar{\rho}}{2}  \|(\bm x-\bm x^*,\bm y-\bm y^*)\|^2 \\
& \mbox{s.t.}\; &  (\bm x,\bm y)\in T_1(\bm x^*,\bm y^*)\cap X(\bm g),
\end{eqnarray*}
 then $(\bm x^*,\bm y^*)$  is an optimal solution  to (CNP) and $\bm x^*$ is an optimal solution to (CNO).

\begin{corollary} %{\bf Corollary 5}
Suppose that $(\bm x^*,\bm y^*)\in X(f)$ and  $f$ is a CN function. Define a set
\begin{eqnarray}
 K_c(\bm x^*,\bm y^*)=\{(\bm x,\bm y) \in R^n\times R^m \mid \nabla g(\bm x^*,\bm y^*)^\top [(\bm x,\bm y)-(\bm x^*,\bm y^*)]<0 \}. \label{d30}
\end{eqnarray} If
%\begin{eqnarray}
$X(\bm g)\cap K_c(\bm x^*,\bm y^*)=\emptyset$ %\label{d31}
%\end{eqnarray}
holds, then $(\bm x^*,\bm y^*)$  is an optimal solution  to (CNP) and $\bm x^*$ is an optimal solution to (CNO).
\end{corollary}

By Corollary 3.2 and Corollary 3.5, if $\bm (\bm x^*,\bm y^*)$  is an optimal solution  to
\begin{eqnarray*}
\mbox{(LCNP)}(\bm x^*,\bm y^*)\qquad & \min\; & \nabla g(\bm x^*,\bm y^*)^\top \bm [(\bm x,\bm y)-(\bm x^*,\bm y^*)]\\
& \mbox{s.t.}\; &  (\bm x,\bm y)\in T_1(\bm x^*,\bm y^*)\cap X(\bm g),
\end{eqnarray*}
 then $(\bm x^*,\bm y^*)$  is an optimal solution  to (CNP) and $\bm x^*$ is an optimal solution to (CNO).

It is clear that $K_w(\bm x^*,\bm y^*)\subset K_c(\bm x^*,\bm y^*)$ and $K_u(\bm x^*,\bm y^*)\subset K_c(\bm x^*,\bm y^*)$.

\begin{theorem}%{\bf Theorem 4}
Suppose that $(\bm x^*,\bm y^*)\in X(f)$ and  $f$ is a weak uniform CN function.  If there is a  neighborhood $O(\bm x^*,\bm y^*)$ of $(\bm x^*,\bm y^*)$   such that
\begin{eqnarray}
X(\bm g)\cap O(\bm x^*,\bm y^*)\cap K_w(\bm x^*,\bm y^*)=\emptyset \label{d32}
\end{eqnarray}
holds,
then $(\bm x^*,\bm y^*)$  is a local optimal solution  to (CNP) and $\bm x^*$ is a local optimal solution to (CNO).
\end{theorem}

%{\bf Proof.} If there is an $(\bm x',\bm y')\in X(\bm g)\cap O(\bm x^*,\bm y^*)$ such that $g(\bm x^*,\bm y^*)> g(\bm x',\bm y')$. We have
%\begin{eqnarray*} 0>g(\bm x',\bm y')-g(\bm x^*,\bm y^*)&\geq& \nabla g(\bm x^*,\bm y^*)^\top\bm [(\bm x',\bm y')-(\bm x^*,\bm y^*)]\\ && +\frac{1}{2} [(\bm x',\bm y')-(\bm x^*,\bm y^*)]^\top B(\bm x^*,\bm y^*) [(\bm x',\bm y')-(\bm x^*,\bm y^*)]. \end{eqnarray*}
%So, $(\bm x',\bm y')\in K_w(\bm x^*,\bm y^*)$. Hence, the conclusion holds.\\

In particular, if $f$ is a strong uniform CN function or $f$ is a uniform CN function,  the conclusion of Theorem 3.4 still holds.\\

\begin{theorem}%{\bf Theorem 5}
Suppose that $(\bm x^*,\bm y^*)\in X(f)$ and  $f$ is a CN function.  If there is a  neighborhood $O(\bm x^*,\bm y^*)$ of $(\bm x^*,\bm y^*)$   such that
\begin{eqnarray}
X(\bm g)\cap O(\bm x^*,\bm y^*)\cap K_c(\bm x^*,\bm y^*)=\emptyset \label{d33}
\end{eqnarray}
holds,
then $(\bm x^*,\bm y^*)$  is a local optimal solution  to (CNP) and $\bm x^*$ is a local optimal solution to (CNO).
\end{theorem}

For $(\bm x^*,\bm y^*)\in X(f)$, let
\begin{eqnarray}
T_0(\bm x^*,\bm y^*)=\{\bm d\in R^n\times R^m \mid \nabla g_i(\bm x^*,\bm y^*)^\top  \bm d= 0,i=1,2,\cdots,r\}.\label{d34}
\end{eqnarray}
For $\bm \alpha \in R^r$, a Lagrange function of (CNP) is defined by
\begin{eqnarray}
L(\bm x,\bm y,\bm \alpha)=g(\bm x,\bm y)+\bm \alpha^\top \bm g(\bm x,\bm y). \label{d35}
\end{eqnarray}

The following necessary conditions of (CNP) are obvious.

\begin{theorem}%{\bf Theorem 6}
Suppose that $(\bm x^*,\bm y^*)\in X(f)$ is a local optimal solution to (CNP). Then

(i) $\nabla g(\bm x^*,\bm y^*)^T\bm d\geq 0$ for $\bm d\in$ $T_0(\bm x^*,\bm y^*)$.

(ii) If $\nabla g_i(\bm x^*,\bm y^*)(i=1,2,\cdots,r)$ is linearly independent, then  there are $\alpha_1^*,\alpha_2^*,\cdots,\alpha_r^* $ such that
\begin{eqnarray}
\nabla g(\bm x^*,\bm y^*)+\sum\limits_{i=1}^r\alpha_i^* \nabla g_i(\bm x^*,\bm y^*)=0. \label{d36}
\end{eqnarray}
Furthermore, if $L(\bm x,\bm y,\bm \alpha^*)$ is convex on $(\bm x,\bm y)$, then $(\bm x^*,\bm y^*)$ is a global optimal solution to (CNP) and $\bm x^*$ is a global optimal solution to (CNO).
\end{theorem}

%{\bf Proof.} (i) The conclusion is clear.

%(ii) By the KKT Theorem, we know that there are $\alpha_1^*,\alpha_2^*,\cdots,\alpha_r^*$ such that the conclusion (\eqref{d36}) holds. If $L(\bm x,\bm y,\bm \alpha^*)$ is convex on $(\bm x,\bm y)$, let any $(\bm x,\bm y)\in X(\bm g)$. So, we have \begin{eqnarray*} L(\bm x,\bm y,\bm \alpha^*)-L(\bm x^*,\bm y^*,\bm \alpha^*)\geq \nabla L(\bm x^*,\bm y^*,\bm \alpha^*)[(\bm x,\bm y)-\bm (x^*,\bm y^*)]. \end{eqnarray*}
%That is $$g(\bm x,\bm y)-g(\bm x^*,\bm y^*)\geq 0.$$
%From the above inequalities,
% $(\bm x^*,\bm y^*)$ is a global optimal solution to (CNP) and $\bm x^*$ is a global optimal solution to (CNO).

In Theorem 3.6, if $\alpha_1^*,\alpha_2^*,\cdots,\alpha_r^*\geq 0$, $(\bm x^*,\bm y^*)$ is a global optimal solution to (CNP) and $\bm x^*$ is a global optimal solution to (CNO).

The following examples  show  that the above sufficient conditions help determine an optimal solution to (CNP) or (CNO), which further help us to design a global algorithm for (CNP), which brings a new way  to study and solve nonconvex and nonsmooth optimization problems.

\begin{example}%{\bf Example 8}
Consider an optimization problem (Example 2.1 in  \cite{Chen1}):
\begin{eqnarray*}
\mbox{(EX8)} &\min& \ f(x_1,x_2)=(x_1+x_2-1)^2+\lambda(|x_1|^{\frac{1}{2}}+|x_2|^{\frac{1}{2}})\\
  &s.t.& \ x_1,x_2\in R^1,
\end{eqnarray*}
where $f(x)$ is a non-smooth and non-convex function. By Example 2.6, $f(x)$ is a weak uniform exact CN  function. Let $(\bm x^*,\bm y^*)\in X(f)$ and $(\bm x,\bm y)=(x_1,x_2,y_1,\cdots,y_6)$.
The weak uniform  exact CN  optimization of (EX8) is defined by
\begin{eqnarray*}
\mbox{(MEX8)}&\min& g(\bm x,\bm y)=(x_1+x_2-1)^2+\lambda(y_1+y_4),\\
&s.t.& (\bm x,\bm y)\in X(f),
\end{eqnarray*}
where $
X(f)=\{(\bm x,\bm y)\mid y_1^4-y_3=0,x_1^2-y_3=0,y_2^2-y_1=0, y_4^4-y_6=0,x_2^2-y_6=0,y_5^2-y_4=0\}$.

%When $\lambda>2$, $(x_1,x_2^*)=(0,0)$ is only an optimal solution to $(EX6)$.
 Then, the linear programming (LCNP-EX8)$(\bm x^*,\bm y^*)$ of (MEX8) at  $(\bm x^*,\bm y^*)$ is  defined by Corollary 3.2. Let $(\bm x^*,\bm y^*)=(0,0,0,0,0,0,0,0)\in X(f)$, $\bm x=(0,0)$ is an optimal solution to (EX8) for $\lambda\geq 8^{\frac{1}{4}}$ in  \cite{Chen1}.
It is clear that there is no optimal solution to (LCNP-EX8) at $(\bm x^*,\bm y^*)=(0,0,0,0,0,0,0,0)\in X(f)$.
Now, a programming (WCNP-EX8)$(\bm x^*,\bm y^*)$ of (MEX8) at  $(\bm x^*,\bm y^*)$ is  defined by Theorem 3.1.
Then, by Theorem 3.3, we compute the objective value of (WCNP-EX8)$(\bm x^*,\bm y^*)$ at $(\bm x^*,\bm y^*)=(0,0,0,0,0,0,0,0)$ such that it is not less than zero. So, we obtian that $\bm x^*=(0,0)$ is an optimal solution to (EX8) for $\lambda\geq \frac{4}{3}\sqrt{\frac{2}{3}}(<8^{\frac{1}{4}})$. Chen et al.(2010)\cite{Chen1} pointed out that there is a smaller error bound $\beta^*<8^{\frac{1}{4}}$ that makes $\bm x^*=(0,0)$ an optimal solution to (EX8). Here, the above  error bound $\beta^*=\frac{4}{3}\sqrt{\frac{2}{3}}$. This result shows that sufficient conditions in Theorem 3 for determining the global optimal solution are valid.

On the other hand, let $\bm d=(d_1,d_2,\cdots,d_8)^\top\in R^8$ for $(\bm x,\bm y)\in R^2\times R^6$. By Theorem 3.1, a programming (WCNP-EX8)$(\bm x^*,\bm y^*)$ of (MEX8) at  $(\bm x^*,\bm y^*)=(0,0,0,0,0,0,0,0)$ is defined.
$(d_1^*,d_2^*,0,0,0,0,0,0)$ is an optimal solution to \mbox{(WCNP-EX8)}$(\bm x^*,\bm y^*)$ when $d_1^*+d_2^*=1$. By (\eqref{d21}) and (\eqref{d22}), we have $\nabla g(\bm x^*,\bm y^*)=(-2,-2,\lambda,0,0,\lambda,0,0)^\top$, $\nabla g_1(\bm x^*,\bm y^*)=(0,0,0,0,-1,0,0,0)^\top$, $\nabla g_2(\bm x^*,\bm y^*)=(0,0,0,0,-1,0,0,0)^\top$, $\nabla g_3(\bm x^*,\bm y^*)=(0,0,-1,0,0,0,0,0)^\top$, $\nabla g_4(\bm x^*,\bm y^*)=(0,0,0,$ $0,0,0,0,-1)^\top$, $\nabla g_5(\bm x^*,\bm y^*)=(0,0,0,0,0,0,0,-1)^\top$, and $\nabla g_6(\bm x^*,\bm y^*)=(0,0,0,0,0,-1,0,0)^\top$. When $(\alpha_1^*,\alpha_2^*,\alpha_3^*,\alpha_4^*,\alpha_5^*,\alpha_6^*)=(0,0,$ $\lambda,0,0,\lambda)$ and $d_1^*+d_2^*=1$ with $\lambda>0$, the optimization condition (\eqref{d21}) and (\eqref{d22}) hold. But, the KKT condition (\eqref{d25}) and (\eqref{d26}) are not true. So, $(\bm x^*,\bm y^*)$ is not a KKT point. $\nabla g_i(\bm x^*,\bm y^*)(i=1,2,\cdots,6)$ is not linear independent.

\end{example}

The above example shows that the optimization condition (\eqref{d21}) and (\eqref{d22}) hold if there is an optimal solution to (CNP) when the objective value of (WCNP-EX8)$(\bm x^*,\bm y^*)$ has a lower bound. The optimization condition of the weak uniform CN form is better that that of CN function in \cite{Jiang}.

\begin{example}%{\bf Example 9.}
Let the function $f(x_1,x_2)=(x_1+x_2-1)^2+\lambda\|(x_1,x_2)\|_0$ be nonconvex and discontinuous, where $\lambda>0$.
\begin{eqnarray*}
\mbox{(EX9)} &\min& \ f(x_1,x_2)=(x_1+x_2-1)^2+\lambda\|(x_1,x_2)\|_0\\
  &s.t.& \ x_1,x_2\in R^1,
\end{eqnarray*}
where $f(x)$ is a non-smooth and non-convex function. $f(x)$ is a weak uniform CN  function. By Example 2.7,
a weak uniform CN  optimization of (EX9) is defined by
\begin{eqnarray*}
\mbox{(EX9)} &\min& \ g(\bm x,\bm y)=(x_1+x_2-1)^2+\lambda(y_1^2+y_2^2),\\
&s.t.&  (\bm x,\bm y)\in X(f),
\end{eqnarray*}
where $(\bm x,\bm y)=(x_1,x_2,y_1,y_2,y_3,y_4)$ and $X(f)$ is defined by Example 2.7. Let $(\bm x^*,\bm y^*)\in X(f)$.
When $(\bm x^*,\bm y^*)=(0,0,0,1,0,1)\in X(f)$, it is clear that $\bm x=(0,0)$ is an optimal solution to (EX9) for $\lambda\geq 2$. So, by Corollary 3.2, it is easily known that there is no optimal solution to (LCNP-EX9)$(\bm x^*,\bm y^*)$ and $X(g)\cap K_c(\bm x^*,\bm y^*)\not=\emptyset$ for $\lambda\geq 1$. But, by Theorem 3.3,
it is clear that for any $(x_1,x_2,y_1,\cdots,y_4)\in X(f)$ and $\lambda\geq 1$ we have
$$-2(x_1+x_2)+\lambda(2y_1+2y_3)+(x_1+x_2)^2\geq -2(x_1+x_2)+1+(x_1+x_2)^2 \geq 0,$$
i.e. $X(f)\cap K_w(\bm x^*,\bm y^*)=\emptyset$. Hence,  $\bm x=(0,0)$ is an optimal solution to (EX9) for $\lambda\geq 1$ by Theorem 3.3. It is easily checked that the optimization condition (\eqref{d21}) and (\eqref{d22}) hold.
\end{example}

Example 3.1 and 3.2 show that the optimization condition (\eqref{d21}) and (\eqref{d22}) probably hold if $f$ a weak uniform CN function, when the KKT condition (\eqref{d25}) and (\eqref{d26}) are not true.

\section{Decomposable Algorithm of (CNP) }

The CN function form $[g:g_1,g_2,\cdots,g_r]$ of $f$ has more variables than $f$.
Because the CN function optimization $\min\limits_{(\bm x,\bm y)\in X(\bm g)} \ g(\bm x,\bm y)$ has more variables than (CNO), it is not easy to solve it, resulting in a scale problem. However, many of decomposable CN functions help reduce  the scale problem of CN function optimization.

Now, we give a decomposable form of  $(\bm x,\bm y))\in X(\bm g)$ of the CN function $f(\bm x)$ on $S$.

Let $f=[g:g_1,g_2,\cdots,g_r]$ be a CN form on $S$. $((\bm x_1,\bm y_1),$
$(\bm x_2,\bm y_2),\cdots,(\bm x_p,\bm y_p))$ is called a decomposition of $(\bm x,\bm y))$ on $S$ if it satisfies the following conditions:

(i) $\bm x=(\bm x_{1}, \bm x_{2},\cdots,\bm x_{p})^\top\in S$, where $\bm x_j=(x_{j1},x_{j2},\cdots,x_{jp_j})\in R^{p_j}, j=1,2,\cdots,p$ and $ \sum\limits_{j=1}^p p_j=n$;

(ii) $\bm y=(\bm y_{1}, \bm y_{2},\cdots,\bm y_{p})^\top\in R^m$, where $\bm y_j=(y_{j1},y_{j2},\cdots,y_{jq_j})\in R^{q_j}, j=1,2,\cdots,p$ and $ \sum\limits_{j=1}^p q_j=m$;

(iii) $((\bm x_1,\bm y_1),(\bm x_2,\bm y_2),\cdots,(\bm x_p,\bm y_p))$ is a rearrangement of $(\bm x,\bm y)$;

(iv) there are no identical variables between $(\bm x_k,\bm y_k)$ and $(\bm x_j,\bm y_j)$ for any $j,k=1,2,\cdots,p, k\not=j)$.

{\bf Note:}  A decomposition of $(\bm x,\bm y)\in X(\bm g)$ means that $(\bm x,\bm y)=((\bm x_1,\bm y_1),(\bm x_2,$
$\bm y_2),\cdots,(\bm x_p,\bm y_p))$. $p$ is called the decomposition number. We have $2\leq p\leq \min\{n,m\}$. $p=2$ is the minimum decomposition number and $p=\min\{n,m\}$ is the maximum decomposition number.

 Let $(\bm x_j,\bm y_j\mid(\bm x,\bm y)):=((\bm x_1,\bm y_1),(\bm x_2,\bm y_2),\cdots,(\bm x_p,\bm y_p))$, where $(\bm x_j,\bm y_j)$ is a variable, i.e. all $(\bm x_k,\bm y_k)(k=1,2,\cdots,p, k\not=j)$ are fixed except  $(\bm x_j,\bm y_j)$.

 %So, $(\bm x_j,\bm y_j)$ is the variable of the function $g(\bm x_j,\bm y_j\mid(\bm x,\bm y)):=g(\bm x,\bm y)$, i.e. all $(\bm x_k,\bm y_k)(k=1,2,\cdots,p, k\not=j)$ are fixed except  $(\bm x_j,\bm y_j)$. $(\bm x_j,\bm y_j)$ is the variable of the function $g_i(\bm x_j,\bm y_j\mid(\bm x,\bm y)):=g(\bm x,\bm y)(i=1,2,\cdots,r)$, i.e. all $(\bm x_k,\bm y_k)(k=1,2,\cdots,p, k\not=j)$ are fixed except  $(\bm x_j,\bm y_j)$.

For each $(\bm x_j,\bm y_j)$, $j=1,2,\cdots,p$, define
\begin{eqnarray*}
X(\bm g_j)=\{(\bm x_j,\bm y_j)\in R^{p_j}\times R^{q_j}\mid  g_i(\bm x_j,\bm y_j\mid(\bm x,\bm y))=0,i=1,2,\cdots,r\}, \label{d7}
\end{eqnarray*}
where $(\bm x,\bm y)\in X(\bm g)$. So, we have
\begin{eqnarray*}
X(\bm g)= X(\bm g_1)\times  X(\bm g_2)\cdots \times X(\bm g_p).\label{d8}
\end{eqnarray*}

\begin{definition} %{\bf Definition 4}
Let the CN  form of $f$ be $[g:g_1,g_2,\cdots,g_r]$ on $S$. If each function $g_1,g_2,\cdots,g_r$ is related only to one of the variables: $\{(\bm x_1,\bm y_1),(\bm x_2,\bm y_2),$ $\cdots,(\bm x_p,\bm y_p)\}$ respectively, i.e.
 the following constraints
\begin{eqnarray*}
g_i(\bm x,\bm y)=g_i((\bm x_1,\bm y_1),(\bm x_2,\bm y_2),\cdots,(\bm x_p,\bm y_p))=0, i=1,2,\cdots,r
\end{eqnarray*}
can be expressed equivalently as
\begin{eqnarray*}
\bm g_j(\bm x_j,\bm y_j)=(g_{j1}(\bm x_{j},\bm y_{j}),g_{j2}(\bm x_{j},\bm y_{j}),\cdots,g_{jr_j}(\bm x_{j},\bm y_{j}))=0, j=1,2,\cdots,p, \label{d13}
\end{eqnarray*}
i.e. \begin{eqnarray} \bm g(\bm x,\bm y)=(\bm g_1(\bm x_1,\bm y_1),\bm g_2(\bm x_2,\bm y_2),\cdots,\bm g_p(\bm x_p,\bm y_p))=0,\label{d13}\end{eqnarray}
then $f$ is called a decomposable CN function on $S$ and $p$ is called the
decomposable number. If there is not any decomposable CN function form, then  $f$ is called an undecomposable CN function. (\eqref{d13}) shows that each function $g_i(\bm x,\bm y)$ is only related to some variable $(\bm x_j,\bm y_j)$ in $\{(\bm x_1,\bm y_1),(\bm x_2,\bm y_2),\cdots,(\bm x_p,\bm y_p)\}$.
\end{definition}

For example, $f$ in Example 2.6 is a decomposable CN function.

For each $j=1,2,\cdots,p$ and $i=1,2,\cdots,r$, let a gradient of $g(\bm x,\bm y)$ and $g_i(\bm x,\bm y)$ on $(\bm x_j,\bm y_j)$ be defined respectively by
$\nabla_j g(\bm x,\bm y):=\nabla_{(\bm x_i,\bm y_j)} g(\bm x_i,\bm y_j\mid(\bm x,\bm y))$ and $\nabla_j g_i(\bm x,\bm y):=\nabla_{(\bm x_i,\bm y_j)} g_i(\bm x_i,\bm y_j\mid(\bm x,\bm y))$,
where all $(\bm x_k,\bm y_k)(k=1,2,\cdots,p, k\not=j)$ are fixed except for $(\bm x_j,\bm y_j)$. We have
\begin{eqnarray}
\nabla g(\bm x,\bm y)=(\nabla_1 g(\bm x,\bm y),\nabla_2 g(\bm x,\bm y),\nabla_p g(\bm x,\bm y))^\top,\label{d14}\\
\nabla g_i(\bm x,\bm y)=(\nabla_1 g_i(\bm x,\bm y),\nabla_2 g_i(\bm x,\bm y),\nabla_p g_i(\bm x,\bm y))^\top. \label{d15}
\end{eqnarray}

 %Let $\bm d\in R^n\times R^m$. Let index set $I=\{1,2,\cdots,r\}$.
%For $(\bm x,\bm y)\in X(\bm g)$, let
%\begin{eqnarray}
%T(\bm x,\bm y)=\{\bm d\in R^n\times R^m \mid \nabla g_i(\bm x,\bm y)^\top  \bm d\leq 0,i=1,2,\cdots,r\}.\label{d16}
%\end{eqnarray}
For $\bm d_j\in R^{p_j}\times R^{q_j},j=1,2,\cdots,p$,   let
\begin{eqnarray}
T_j(\bm x_j,\bm y_j)=\{ \bm d_j\in R^{p_j}\times R^{q_j}\mid \nabla_j g_{i}(\bm x,\bm y)^\top \bm d_j\leq 0,i=1,2\cdots,r\}. \label{d17}
\end{eqnarray}
By (\eqref{d17}), we have
\begin{eqnarray*}
T_1(\bm x_1,\bm y_1)\times T_j(\bm x_2,\bm y_2)\cdots \times T_p(\bm x_p,\bm y_p) \subset T(\bm x,\bm y). \label{d14}
\end{eqnarray*}

By Definition 4.1, we have following propositions.

\begin{proposition}%{\bf Proposition 14}
Let $f$ be a decomposable CN function and $f=[g:g_1,g_2,\cdots,$ $g_r]$,
where $((\bm x_1,\bm y_1),(\bm x_2,\bm y_2),\cdots,(\bm x_p,\bm y_p))$ is a decomposition of $(\bm x,\bm y)$.
Then
\begin{eqnarray}
T_1(\bm x_1,\bm y_1)\times T_2(\bm x_2,\bm y_2)\cdots \times T_p(\bm x_p,\bm y_p)= T(\bm x,\bm y). \label{d18}
\end{eqnarray}

%{\bf Proof.} By Definition 4 and $\bm d=(\bm d_1,\bm d_2,\cdots, \bm d_p)\in T(\bm x,\bm y)$, we have $\bm d_j \in T_j(\bm x_j,\bm y_j),j=1,2,\cdots,p$. Hence, (\eqref{d18}) is true.

\end{proposition}

\begin{proposition}%{\bf Proposition 15}
 Let  $f=[g:g_1,g_2,\cdots,g_r]$ be a decomposable CN function and
 $((\bm x_1,\bm y_1),(\bm x_2,\bm y_2),\cdots,(\bm x_p,\bm y_p))$ be a decomposition of $(\bm x,\bm y)$.
If $g(\bm x,\bm y)$ is a (weak,strong) uniform convex function on $(\bm x,\bm y)$, then $g(\bm x_i,\bm y_j\mid(\bm x,\bm y))$ is a (weak,strong) uniform convex function on $(\bm x_j,\bm y_j)(j=1,2,\cdots,p)$.%, where $(\bm x_j,\bm y_j)$ is variable. %i.e. all $(\bm x_k,\bm y_k)(k=1,2,\cdots,p, k\not=j)$ are fixed except for $(\bm x_j,\bm y_j)$.

%(ii) If $g(\bm x,\bm y)$ is a strong uniform convex function on $(\bm x,\bm y)$, then $g(\bm x_i,\bm y_j\mid(\bm x,\bm y))$ is a strong uniform convex function on $(\bm x_j,\bm y_j)(j=1,2,\cdots,p)$.%, where $(\bm x_j,\bm y_j)$ is variable. %i.e. all $(\bm x_k,\bm y_k)(k=1,2,\cdots,p, j\not=k)$ are fixed except for $(\bm x_j,\bm y_j)$.

%(iii) $g(\bm x,\bm y)$ is an uniform convex function on $(\bm x,\bm y)$ if and only if $g(\bm x_i,\bm y_j\mid(\bm x,\bm y))$ is an uniform convex function on $(\bm x_j,\bm y_j)(j=1,2,\cdots,p)$. %, where $(\bm x_j,\bm y_j)$ is variable.
 %i.e. all $(\bm x_k,\bm y_k)(k=1,2,\cdots,p, k\not=j)$ are fixed except for $(\bm x_j,\bm y_j)$.

%(iv) $g(\bm x,\bm y)$ is a convex function on $(\bm x,\bm y)$ if and only if $g(\bm x_i,\bm y_j|(\bm x,\bm y))$ is a convex function on $(\bm x_j,\bm y_j)(j=1,2,\cdots,p)$.%, where $(\bm x_j,\bm y_j)$ is variable.
 %i.e. all $(\bm x_k,\bm y_k)(k=1,2,\cdots,p, k\not=j)$ are fixed except for $(\bm x_j,\bm y_j)$.
\end{proposition}

An example of decomposable CN function is given as follows.

\begin{example}
(In \cite{Chen2}) The function in sparse optimization is
\begin{eqnarray}
f(\bm x)=\|A\bm x-\bm b\|^2+\lambda\|\bm x\|_0,\label{d10}
\end{eqnarray}
where $A\in R^m\times R^n$, $\bm b\in R^m$, $\lambda>0$ and $\|\bm x\|_0$ is 0-norm. Then, a weak uniform  CN function form of $f(\bm x)$  is obtained by
\begin{eqnarray*}
g(\bm x,\bm y)&=&\|A\bm x-\bm b\|^2+\lambda\sum\limits_{i=1}^n y_i^2 :\\
g_i(\bm x,\bm y)&=& (x_i+y_i-1)^2-y_{i+n}=0,\ \ i=1,2,\cdots, n,\\
g_{i+n}(\bm x,\bm y)&=& x_i^2+(y_i-1)^2-y_{i+n}=0,\ \ i=1,2,\cdots, n, \\
g_{i+2n}(\bm x,\bm y)&=&y_i^2-y_i=0,\ \ i=1,2,\cdots, n,
\end{eqnarray*}
where $\bm y\in R^{2n}$. Let $\bm w_i=(x_i,y_i,y_{i+n})$, $i=1,2,\cdots,n$. Then, a decomposable CN form of $f(\bm x)$  at the maximum  decomposition number $p=n$  is defined by
\begin{eqnarray*}
g(\bm w_1,\bm w_2,\cdots,\bm w_n)&=&\|A\bm x-\bm b\|^2+\lambda\sum\limits_{i=1}^n y_i^2:\\
\bm g_i(\bm w_i)&=&((x_i+y_i-1)^2-y_{i+n},x_i^2+(y_i-1)^2-y_{i+n},y_i^2-y_i)\\
&=&\bm 0,\  i=1,2,\cdots, n.
\end{eqnarray*}
\end{example}

If $((\bm x_1,\bm y_1),$
$(\bm x_2,\bm y_2),\cdots,(\bm x_p,\bm y_p))$ is a decomposition of $(\bm x,\bm y))$, (CNP) is redefined as
\begin{eqnarray*}
\mbox{(CNP)}\qquad & \min\; & g(\bm x,\bm y)=g(\bm x_1,\bm y_1),(\bm x_2,\bm y_2),\cdots,(\bm x_p,\bm y_p)) \\
& \mbox{s.t.}\; & \bm g_j(\bm x_j,\bm y_j|(\bm x,\bm y))=0,j=1,2,\cdots,p,\\
&& (\bm x_j,\bm y_j)\in R^{p_j}\times R^{q_j},j=1,2,\cdots,p,
\end{eqnarray*}
where $f$ is not necessarily a decomposable CN function. For each
$j=1,2,\cdots,p$, the $j$th subproblem of (CNP) is defined by
\begin{eqnarray*}
\mbox{(CNP)}_j\qquad & \min\; & g(\bm x_i,\bm y_j|(\bm x,\bm y)) \\
& \mbox{s.t.}\; & (\bm x_i,\bm y_j|(\bm x,\bm y))\in X(\bm g_j),
\end{eqnarray*}
where $(\bm x_j,\bm y_j)$ is the variable, i.e. all $(\bm x_k,\bm y_k)(k=1,2,\cdots,p, k\not=j)$ are fixed except for $(\bm x_j,\bm y_j)$ in problem (CNP)$_j$. Then, the optimal solution to (CNP) is expected to be obtained by solving $p$ subproblems (CNP)$_j$, $j=1,2,\cdots,p.$ If $f$ is a decomposable CN function and there is an optimal solution to (CNP) or (CNO), there are optimal solutions to all subproblems  (CNP)$_1$,(CNP)$_2$,$\cdots$, (CNP)$_p$ of (CNP). But, if there is not an optimal solution to anyone of all subproblems  (CNP)$_1$,(CNP)$_2$,$\cdots$, (CNP)$_p$ of (CNP), there is not an optimal solution to (CNP) or (CNO).\\

The following theorems are true.

\begin{theorem}%{\bf Theorem 7}
Suppose that $(\bm x^*,\bm y^*)=((\bm x_1^*,\bm y_1^*),(\bm x_2^*,\bm y_2^*),\cdots,(\bm x_p^*,\bm y_p^*))\in X(f)$ and  $f$ is a decomposable CN function. For each $j=1,2,\cdots,p$,
let the problem
\begin{eqnarray*}
\mbox{(WCNP)}_j(\bm x^*_j,\bm y^*_j)\qquad & \min\; & \nabla_j g(\bm x^*,\bm y^*)^\top \bm d_j \\
& \mbox{s.t.}\; &  \bm d_j\in T_j(\bm x^*_j,\bm y^*_j).
\end{eqnarray*}
For all $j=1,2,\cdots,p$,  if $\bm d^*_j$  is an optimal solution  to (WCNP)$_j(\bm x^*_j,\bm y^*_j)$ such that $\nabla_j g(\bm x^*,\bm y^*)^\top \bm d^*_j  \geq 0$, then $(\bm x^*,\bm y^*)$  is an optimal solution  to (CNP) and $\bm x^*$ is an optimal solution to (CNO). Furthermore, if there is a $\bm d'$ such that $\nabla g_i(\bm x^*,\bm y^*)^\top  \bm d'< 0,i=1,2,\cdots,r,$ then there are $\alpha_i^*\geq 0, i=1,2,\cdots,r$ such that
  \begin{eqnarray}
 \nabla g(\bm x^*,\bm y^*) + \sum\limits_{i=1}^r \alpha_i^* \nabla g_i(\bm x^*,\bm y^*)=0. \label{d37}
\end{eqnarray}
\end{theorem}

%{\bf Proof.} Since $ T(\bm x^*,\bm y^*)=\bigcup\limits_{j=1}^p T_j(\bm x^*_j,\bm y^*_j)$, we have $\bm d^*=(\bm d_1^*,\bm d_2^*,\cdots,\bm d_p^*)\in T(\bm x^*,\bm y^*)$. So,
%$$\nabla g(\bm x^*,\bm y^*)^\top \bm d^*=\sum\limits_{j=1}^p \nabla_j g(\bm x^*,\bm y^*)^\top \bm d_j^*\geq 0.$$ By Corollary 2, the conclusion is true.

Theorem 4.1 shows that the optimal solution to (CNP) is expected to be obtained by solving $p$ subproblems (CNP)$_j$, $j=1,2,\cdots,p.$ These subproblems (CNP)$_j(j=1,2,\cdots,p)$ have smaller scale although the problem (CNP) has more variables.

\begin{theorem}%{\bf Theorem 8}
Suppose that $(\bm x^*,\bm y^*)=((\bm x_1^*,\bm y_1^*),(\bm x_2^*,\bm y_2^*),\cdots,(\bm x_p^*,\bm y_p^*))\in X(f)$ and  $f$ is a decomposable uniform CN function. For each $j=1,2,\cdots,p$,
let the problem
\begin{eqnarray*}
\mbox{(UCNP)}_j(\bm x^*_j,\bm y^*_j)\qquad & \min\; & \nabla_j g(\bm x^*,\bm y^*)^\top \bm d_j+\frac{\bar{\rho}}{2}\|\bm d_j\|^2 \\
& \mbox{s.t.}\; &  \bm d_j\in T_j(\bm x^*_j,\bm y^*_j).
\end{eqnarray*}
For all $j=1,2,\cdots,p$,  if $\bm d^*_j$  is an optimal solution  to (UCNP)$_j(\bm x^*_j,\bm y^*_j)$ such that $\nabla_j g(\bm x^*,\bm y^*)^\top \bm d^*_j+\frac{\bar{\rho}}{2}\|\bm d_j^*\|^2   \geq 0$, then $(\bm x^*,\bm y^*)$  is an optimal solution  to (CNP) and $\bm x^*$ is an optimal solution to (CNO). Furthermore, if there is a $\bm d'$ such that $\nabla g_i(\bm x^*,\bm y^*)^\top  \bm d'< 0,i=1,2,\cdots,r,$ then there are $\alpha_i^*\geq 0, i=1,2,\cdots,r$ such that
  \begin{eqnarray}
 \frac{\bar{\rho}}{2}\bm d_j^*+\nabla g(\bm x^*,\bm y^*) + \sum\limits_{i=1}^r \alpha_i^* \nabla g_i(\bm x^*,\bm y^*)=0. \label{d38}
\end{eqnarray}
\end{theorem}

\begin{theorem}%{\bf \bf Theorem 9}
 Suppose that $(\bm x^*,\bm y^*)\in X(f)$ and  $f$ is an uniform decomposable CN function. For each $j=1,2,\cdots,p$, define a set
\begin{eqnarray}
 K_{uj}(\bm x^*_j,\bm y^*_j)=\{(\bm x_j,\bm y_j) \in R^{p_j}\times R^{q_j} \mid \nabla_j g(\bm x^*,\bm y^*)^\top (\bm x_j-\bm x^*_j,\bm y_j-\bm y^*_j)\nonumber\\
 +\frac{\bar{\rho}}{2} \|(\bm x_j-\bm x^*_j,\bm y_j-\bm y^*_j)\|^2<0 \}. \label{d39}
\end{eqnarray} If
%\begin{eqnarray}
$X(\bm g_j)\cap K_{uj}(\bm x^*_j,\bm y^*_j)=\emptyset$ %\label{d40}
%\end{eqnarray}
 holds, then $(\bm x^*,\bm y^*)$  is an optimal solution  to (CNP) and $\bm x^*$ is an optimal solution to (CNO).
\end{theorem}
%{\bf Proof.} By (\eqref{d39}) and (\eqref{d40}), we have $K_{u}(\bm x^*,\bm y^*)=\bigcup\limits_{j=1}^p K_{uj}(\bm x^*_j,\bm y^*_j)$. So, $X(\bm g)\cap K_{u}(\bm x^*,\bm y^*)=\emptyset$. By Corollary 4,  $(\bm x^*,\bm y^*)$  is an optimal solution  to (CNP) and $\bm x^*$ is an optimal solution to (CNO).\\

%Similarly, the following conclusion holds.
\begin{theorem} %{\bf Theorem 10}
Suppose that $(\bm x^*,\bm y^*)\in X(f)$ and  $f$ is a decomposable CN function.  For each $j=1,2,\cdots,p$, define a set
\begin{eqnarray}
{\small K_{cj}(\bm x^*_j,\bm y^*_j)=\{(\bm x_j,\bm y_j) \in R^{p_j}\times R^{q_j} \mid \nabla_j g(\bm x^*,\bm y^*)^\top (\bm x_j-\bm x^*_j,\bm y_j-\bm y^*_j)<0 \}.} \label{d41}
\end{eqnarray} If
%\begin{eqnarray}
$X(\bm g_j)\cap K_{cj}(\bm x^*_j,\bm y^*_j)=\emptyset $ % \label{d42}
%\end{eqnarray}
 holds, then $(\bm x^*,\bm y^*)$  is an optimal solution  to (CNP) and $\bm x^*$ is an optimal solution to (CNO).\\
\end{theorem}

Now, an augmented Lagrange penalty function for (CNP) with decomposable variable is defined. Next, an algorithm is proposed by the augmented Lagrange penalty function for (CNP) and their convergence is proved.

Let $\bm \alpha_j\in R^{r_j},j=1,2,\cdots,p$ be Lagrange parameters and $\sigma>0$ be a penalty parameter  and $\bm \alpha=(\bm \alpha_1,\bm \alpha_2,\cdots,\bm \alpha_p)$. Suppose that $((\bm x_1,\bm y_1),(\bm x_2,\bm y_2),\cdots,(\bm x_p,\bm y_p))$ is a decomposable variable of $(\bm x,\bm y)$ and $f$  is not necessarily a decomposable CN function.
In order to solve (CNP), the augmented Lagrange penalty functions for all subproblems (CNP)$_j(j=1,2,\cdots,p)$ are defined by
\begin{eqnarray}
 A_j(\bm x_i,\bm y_j\mid(\bm x,\bm y);\bm \alpha_j,\sigma)&=&g(\bm x_i,\bm y_j|(\bm x,\bm y))+\bm \alpha_j^\top \bm g_{j}(\bm x_{j},\bm y_{j}|(\bm x,\bm y))\nonumber \\
&& +\frac{1}{2}\sigma \|\bm g_{j}(\bm x_{j},\bm y_{j}|(\bm x,\bm y))\|^2, \label{d43}
\end{eqnarray}
where  $(\bm x_j,\bm y_j)$ is variable, i.e. all $(\bm x_k,\bm y_k)(k=1,2,\cdots,p, k\not=j)$ are fixed except for $(\bm x_j,\bm y_j)$. By (\eqref{d43}), for $j=1,2,\cdots,p$, define an unconstraint optimization problem
\begin{eqnarray*}
\mbox{(CNP)}_j(\bm \alpha_j,\sigma)\qquad & \min\; & A_j(\bm x_j,\bm y_j|(\bm x,\bm y);\bm \alpha_j,\sigma) \\
& \mbox{s.t.}\; & (\bm x_j,\bm y_j)\in R^{p_j}\times R^{q_j}.
\end{eqnarray*}

To solve the problem (CNP)$_j(\bm \alpha_j,\sigma)$, an algorithm involving an augmented Lagrange penalty function for (CNP)(which is called Algorithm 1) is proposed.\\

{\bf Algorithm 1:}

\begin{description}
\item[Step 1:] Let $\epsilon>0,\sigma_{1}>0,N>1, \bm w^0_j=(\bm x_j^0,\bm y_j^0)\in R^{p_j}\times R^{q_j}(j=1,2,\cdots,p), \bm \alpha_j^1\in R^{r_j}(j=1,2,\cdots,p)$, $k=1$.

\item[Step 2.1:] Let $j=1$. If a point
     $$(\bm x^{k-1},\bm y^{k-1})_1=((\bm x_1^{k-1},\bm y_1^{k-1}),(\bm x_2^{k-1},\bm y_2^{k-1}),\cdots,(\bm x_{p}^{k-1},\bm y_{p}^{k-1}))$$
    is obtained, then find $(\bm x_1^k,\bm y_1^k)\in R^{p_1}\times R^{q_1} $ to the subproblem $\min\limits_{(\bm x_1,\bm y_1)}A_1^k(\bm x_1,$ $\bm y_1|(\bm x^{k-1},\bm y^{k-1})_1;\bm \alpha_1^k,\sigma_k)$ such that
    $$\nabla_1 A_1^k(\bm x_1^{k},\bm y_1^{k}|(\bm x^{k-1},\bm y^{k-1})_1;\bm \alpha_1^{k},\sigma_{k})=0,\ \ $$ where $(\bm x_1,\bm y_1)$ of function $A_1^k(\bm x_1,\bm y_1|(\bm x^{k-1},\bm y^{k-1})_1;\bm \alpha_1^k,\sigma_k)$ is  variable, i.e. all $(\bm x_s^{k-1},\bm y_s^{k-1})(s=2,\cdots,p)$ are fixed except for $(\bm x_1,\bm y_1)$. Let $j=2$ and go to Step 2.2.

\item[Step 2.2:] Let $j>1$. If a point
     $$(\bm x^{k-1},\bm y^{k-1})_j=((\bm x_1^{k},\bm y_1^{k}),\cdots,(\bm x_{j-1}^{k},\bm y_{j-1}^{k}),(\bm x_{j}^{k-1},\bm y_{j}^{k-1}),\cdots,(\bm x_{p}^{k-1},\bm y_{p}^{k-1}))$$
    is obtained, then find $(\bm x_{j}^k,\bm y_{j}^k)\in R^{p_j}\times R^{q_j} $ to the subproblem $\min\limits_{(\bm x_j,\bm y_j)}A_j^k(\bm x_j,$ $\bm y_j|(\bm x^{k-1},\bm y^{k-1})_j;\bm \alpha_j^k,\sigma_k)$ such that
    $$\nabla_j A_j^k(\bm x_j^{k},\bm y_j^{k}|(\bm x^{k-1},\bm y^{k-1})_j;\bm \alpha_j^{k},\sigma_{k})=0,\ \ $$ where
 $(\bm x_j,\bm y_j|(\bm x^{k-1},\bm y^{k-1})_j=((\bm x_1^{k},\bm y_1^{k}),\cdots,(\bm x_{j-1}^{k},\bm y_{j-1}^{k}),(\bm x_{j},\bm y_{j}),(\bm x_{j+1}^{k-1},$ $\bm y_{j+1}^{k-1}),\cdots,(\bm x_{p}^{k-1},\bm y_{p}^{k-1})),$ i.e., $(\bm x_j,\bm y_j)$ of $A_j^k(\bm x_j,\bm y_j|(\bm x^{k-1},\bm y^{k-1})_j;\bm \alpha_j^k,\sigma_k)$ is  variable, i.e. all $(\bm x_s^{k},\bm y_s^{k})(s=1,2,\cdots,j-1)$ and $(\bm x_s^{k-1},\bm y_s^{k-1})(s=j+1,j+2,\cdots,p)$ are fixed except for $(\bm x_j,\bm y_j)$. Go to Step 2.3

    \item[Step 2.3:]  If $j=p$ and go to Step 3.  Otherwise, $j:=j+1$ and go to Step 2.2.

\item[Step 3:]If $(\bm x^k,\bm y^k)=(\bm x^{k-1},\bm y^{k-1})\in X(f)$ and  $L(\bm x,\bm y,\bm \alpha^{k})$ is convex on $(\bm x,\bm y)$ (see (16)), then stop and  $\bm x^k$ is an optimal solution to (CNO). Otherwise, go to Step 4.

\item[Step 4:] If $\|\bm g(\bm x^k,\bm y^k)\|=\sum_{j=1}^p\|\bm g_{j}(\bm x_{j}^k,\bm y_{j}^k)\|<\epsilon$, then stop and $\bm x^k$ is an approximate solution to (CNO). Otherwise, for $j=1,2,\cdots,p$,
     let $\bm \alpha_j^{k+1}=\bm \alpha_j^{k}+\sigma_{k}\bm g_j(\bm x_j^{k},\bm y_j^{k})$, $\sigma_{k+1}=N\sigma_{k}$, $k:=k+1$ and go to Step 2.1.\\
\end{description}

{\bf Note:} In Step 3, if $(\bm x^k,\bm y^k)=(\bm x^{k-1},\bm y^{k-1})\in X(f)$ holds,  we have $(\bm x_{j}^k,\bm y_{j}^k)\in R^{p_j}\times R^{q_j}(j=1,2,\cdots,p)$. So, for all $(j=1,2,\cdots,p$, $$\nabla_j A_j^k(\bm x_j^{k},\bm y_j^{k}|(\bm x^{k-1},\bm y^{k-1})_j;\bm \alpha_j^{k},\sigma_{k})=\nabla_j A_j^k(\bm x_j^{k},\bm y_j^{k}|(\bm x^{k},\bm y^{k})_j;\bm \alpha_j^{k},\sigma_{k})=0$$ holds.

By Algorithm 1, if $\|\bm g(\bm x^k,\bm y^k)\|<\epsilon$ holds, it maybe to be able to find an approximate global optimal solution to (CNO). Under some conditions, it is proved that Algorithm 1 can converge to a KKT point for $\epsilon=0$.

Let
$$S(\pi, g)=\{(\bm x,\bm y)\mid \pi\geq  g(\bm x,\bm y)\}, $$
which is called a level set.  If $S(\pi, g)$ is bounded for any given $\pi>0$, then $S(\pi, g)$ is also bounded.

\begin{theorem}
%{\bf Theorem 11}
 Let $\epsilon=0$ and $f$ be a decomposable CN function. Suppose that a sequence of $\{(\bm x^k,\bm y^k):=((\bm x_1^k,\bm y_1^k),(\bm x_2^k,$ $\bm y_2^k),\cdots,(\bm x_p^k,\bm y_p^k))\}$, $k=1,2,\cdots$, is obtained by Algorithm 1.
Let the sequence of $\{H_k(\bm x^{k},\bm y^{k},\sigma_{k})\}$, $k=1,2,\cdots,$ be bounded and the level set $S(\pi, g)$ be bounded, where
$$H_k(\bm x^{k},\bm y^{k},\rho_{k})=g(\bm x^{k},\bm y^{k})+\sigma_{k}\sum\limits_{i=1}^r g_i(\bm x^{k},\bm y^{k})^2.$$
(i) If the algorithm stops at a finite number of step $k$, then $\bm x^k$  is a global optimal solution to (CNO).\\
(ii) If the sequence $\{(\bm x^{k},\bm y^{k})\}$ is an infinite sequence, then $\{(\bm x^{k},\bm y^{k})\}$ is bounded and any limit point $(\bm x^*,{\bm y}^*)$ of the sequence belongs to $X(\bm g)$, and
there exist $\eta>0$ and $\lambda_i$, $i=1, 2, \cdots, r$, such that
\begin{eqnarray}
\eta\nabla g(\bm x^*,{\bm y}^*)+ \sum\limits_{i=1}^r \lambda_i \nabla g_i(\bm x^*,{\bm y}^*)=0. \label{d44}
\end{eqnarray}
 If  $(\bm x^*,{\bm y}^*)\in X(f)$ and $\eta g(\bm x,{\bm y})+ \sum\limits_{i=1}^r \lambda_i g_i(\bm x,{\bm y})$ is convex on $(\bm x,\bm y)$ or $\lambda_i\geq 0$, $i=1, 2, \cdots, r$, then
$\bm x^*$ is an optimal solution to (CNO).
\end{theorem}

{\it Proof.} (i) By Step 2 and Step 3 of Algorithm 1, we have
\begin{eqnarray*}
 \nabla_j A_j^k(\bm x_j^{k},\bm y_j^{k}|(\bm x^{k-1},\bm y^{k-1})_j;\bm \alpha_j^{k},\sigma_{k})&=&\nabla_j A_j^k(\bm x_j^{k},\bm y_j^{k}|(\bm x^{k},\bm y^{k})_j;\bm \alpha_j^{k},\sigma_{k})\\
 &=&\nabla_j g_j(\bm x^k,{\bm y}^k)+ {\bm \alpha_j^{k}}^\top \nabla_j \bm g_j(\bm x_j^k,{\bm y_j}^k)\\
 &=& 0,\ \ \ \ j=1,2,\cdots,p,
 \end{eqnarray*}
where $\bm \alpha^{k}=(\bm \alpha_1^{k},\bm \alpha_2^{k},\cdots, \bm \alpha_p^{k})^\top$.
So,
%\begin{eqnarray*}
$\nabla g(\bm x^k,{\bm y}^k)+ {\bm \alpha^{k}}^\top \nabla \bm g(\bm x^k,{\bm y}^k)=0.$
% \end{eqnarray*}
By Theorem 4.1, $\bm x^k$  is a global optimal solution to (CNO).

(ii) By Algorithm 1,  as $k\to +\infty$, since $\{H_k(\bm x^k,{\bm y}^k,\rho_k)\}$ is bounded,  there must be some $\pi>0$ such
that\begin{eqnarray*}
          \pi>H_k(\bm x^k,{\bm y}^k,\sigma_k)\geq  g(\bm x^k,{\bm y}^k).
 %         &=&  g(\bm x^k,{\bm y}^k)+\sigma_k\sum\limits_{i=1}^r g_i(\bm x^k,{\bm y}^k)^2\\
 %         &\geq&  g(\bm x^k,{\bm y}^k).
\end{eqnarray*}
 $\{(\bm x^k,{\bm y}^k)\}$ is bounded because the level set $S(\pi,f)$ is bounded.
Without loss of generality, suppose $(\bm x^k,{\bm y}^k)\to (\bm x^*,{\bm y}^*)$, i.e. $(\bm x_j^k,\bm y_j^k)\to (\bm x_j^*,\bm y_j^*) (j=1,2,\cdots,p)$.  Since $g$ is continuous, $S(\pi,g)$ is closed. So, $g(\bm x^k,{\bm y}^k)$ is bounded and there is a $\sigma'>0$ such that $g(\bm x^k,{\bm y}^k)> -\sigma'$.

From the above inequality, we have
that\begin{eqnarray*}
       \sum\limits_{i=1}^r g_i(\bm x^k,{\bm y}^k)^2 \leq\frac{1}{\sigma_k} (\pi- g(\bm x^k,{\bm y}^k))<\frac{\pi+\sigma'}{\sigma_k}.
\end{eqnarray*}
We have $\sum\limits_{i=1}^r (g_i(\bm x^k,{\bm y}^k))^2\to 0$ as  $\sigma_k\to +\infty$.
So, $(\bm x^*,{\bm y}^*)\in X(\bm g)$.

By  Algorithm 1,  there is an infinite sequence $\{(\bm x^k,{\bm y}^k,\bm \alpha^k,\rho_k) \}$ such that
$$\nabla_j A_j^k(\bm x_j^{k},\bm y_j^{k}|(\bm x^{k-1},\bm y^{k-1})_j;\bm \alpha_j^{k},\sigma_{k})=0,\  j=1,2,\cdots,p.$$
For all $j=1,2,\cdots,p$, we have
 \begin{eqnarray}
\nabla_j g(\bm x_j^{k},\bm y_j^{k}|(\bm x^{k-1},\bm y^{k-1})_j)
+ \sum\limits_{i=1}^{r_j} \alpha_{ji}^{k+1}\nabla_j g_{ji}(\bm x_j^{k},\bm y_j^{k})=0,\ \label{d45}
\end{eqnarray}
where $\bm \alpha_j^{k}=(\alpha_{j1}^{k},\alpha_{j2}^{k},\cdots,\alpha_{jr_j}^{k})$,
$(\bm x_j^{k},\bm y_j^{k}|(\bm x^{k-1},\bm y^{k-1})_j)=((\bm x_1^{k},\bm y_1^{k}),\cdots,(\bm x_{j-1}^{k},$ $\bm y_{j-1}^{k}),(\bm x_{j}^k,\bm y_{j}^k),(\bm x_{j+1}^{k-1},$ $\bm y_{j+1}^{k-1}),\cdots,(\bm x_{p}^{k-1},\bm y_{p}^{k-1}))$
and $\alpha_{ji}^{k+1}=\alpha_{ji}^{k}+\sigma_{k} g_{ji}(\bm x_{ji}^{k},\bm y_{ji}^{k})$ $(i=1,2,\cdots,r_j)$.
For all $j=1,2,\cdots,p$, let
 \begin{eqnarray*}
\gamma_j^k=1+\sum\limits_{i=1}^{r_j} (\max\{\alpha_{ji}^{k+1},0\}+\max\{-\alpha_{ji}^{k+1},0\})>0.
\end{eqnarray*}
 Let $\eta_j^k=\frac{1}{\gamma_j^k}>0$,
$\mu_{ji}^k=\frac{\max\{\alpha^{k+1}_{ji},0\}}{\gamma_j^k}\geq 0, i=1,2,\cdots,r_j$ and $\nu_{ji}^k=\frac{\max\{-\alpha_{ji}^{k+1},0\}}{\gamma_j^k}\geq 0,i=1,2,\cdots,r_j$.
Then,
 \begin{eqnarray}
 \eta_j^k+\sum\limits_{i=1}^{r_j}(\mu_{ji}^k+\nu_{ji}^k)=1.\label{d46}
 \end{eqnarray}
Clearly, for all $j=1,2,\cdots,p$, as $k\to\infty$, we have $\eta_j^k\rightarrow \eta_j>0, \mu_{ji}^k\rightarrow \mu_{ji}, \nu_{ji}^k\rightarrow \nu_{ji},\forall i=1,2,\cdots,r_j$.
By (\eqref{d45}) and (\eqref{d46}), we have
 \begin{eqnarray}
\eta\nabla g(\bm x^*,{\bm y}^*)+ \sum\limits_{i=1}^r( \mu_i -\nu_i)\nabla g_i(\bm x^*,{\bm y}^*)=0.\label{d47}
\end{eqnarray}
By (\eqref{d47}), let $\lambda^k=\mu_i^k-\nu_i^k\to \lambda$ as $k\rightarrow +\infty$,  and we have (\eqref{d44}).\\

 Finally, four examples are given to show Algorithm 1 can solve an approximate optimal solution to (CNO). All codes are written with Matlab2016a and numerical experiments are carried out with Thinkpad S3.
 % by Matlab2016a.

%\end{example}
\begin{example}
%{\bf Example 11.}
Consider the optimization problem which is the same as Example 3.1. $(x_1^*,x_2^*)=(0,0)$ is only an optimal solution to (Ex4.2) and $f(0,0)=1$ for $\lambda>2$. The CN function $f(x)$ is decomposable.
By Example 3.1 and Algorithm 1, two subproblems are solved:
\begin{eqnarray*}
\mbox{(Ex4.2)}_1 &\min& \ g((x_1,y_1,y_2,y_3),(x_2,y_4,y_5,y_6))=(x_1+x_2-1)^2+\lambda(y_2^2+y_5^2)\\
  &s.t.& \bm g_1(x_1,\bm y_1)=(y_1^4-y_3,x_1^2-y_3,y_2^2-y_1)=0.
\end{eqnarray*}
\begin{eqnarray*}
\mbox{(Ex4.2)}_2 &\min& \ g((x_1,y_1,y_2,y_3),(x_2,y_4,y_5,y_6))=(x_1+x_2-1)^2+\lambda(y_2^2+y_5^2)\\
  &s.t.& \bm g_2(x_2,\bm y_2)=(y_4^4-y_6,x_2^2-y_6,y_5^2-y_4))=0.
\end{eqnarray*}
Now, let stating parameters $\lambda=20000,\epsilon=10^{-4}, \theta=0.1,\sigma_1=1000, N=1000, \alpha_1^1=(2,2,2)^\top,\alpha_1^2=(2,2,2)^\top$ and $((x_1^0,y_1^0,$ $y_2^0,y_3^0),(x_2^0,y_4^0,y_5^0,y_6^0))=((2,2,2,2),(2,2,2,2))$. At iteration 3, an approximate optimal solution is obtained $((x_1^k,y_1^k,y_2^k,y_3^k),$ $(x_2^k,y_4^k,y_5^k,y_6^k))=((0.0009,0.0000,0.0000, 0.0000), $\\ $ (0.0006,0.0000,-0.0000,0.0000))$.  Furthermore, when  random  values of starting points are chosen from $[-500,500]$, the same approximate optimal solution is obtained at iteration 6.
\end{example}

\begin{example}
Consider the special sparse optimization problem in \cite{Chen2}:
 \begin{eqnarray*}
\mbox{(Ex4.3)} &\min& f_n(\bm x)=(\sum\limits_{i=1}^n ix_i-2n)^2+\lambda\sum\limits_{i=1}^n |x_i|_0 \\
&s.t.& \bm x\in R^n.
\end{eqnarray*}
 For $(x_i,y_i,y_{i+n})$, $i=1,2,\cdots,n$, a decomposable PCN form of $f_n(\bm x)$  at the maximum  decomposition number $p=n$  is defined by
\begin{eqnarray*}
[(\sum\limits_{i=1}^n ix_i-2n)^2+\lambda\sum\limits_{i=1}^n y_{i}^2 :
((x_i+y_i-1)^2-y_{i+n},\\
x_i^2+(y_i-1)^2-y_{i+n},y_{i}^2-y_i)=0,  i=1,2,\cdots,n].
\end{eqnarray*}
Let $n=ep$,  $\bm w_{j}=(x_{ej+1},x_{ej+2},\cdots,x_{ej+e},y_{ej+1},y_{ej+2},\cdots, y_{ej+e},y_{ej+1+n},$ $ y_{ej+2+n},\cdots,y_{ej+e+n})$,
$j=0,1,2,\cdots,p-1$. So,for $j=0,1,2,\cdots,p-1$, $p$ subproblems are solved by:
 \begin{eqnarray*}
\mbox{(Ex4.3)}_j &\min& g(\bm w_0,\bm w_1,\cdots,\bm w_{p-1})=(\sum\limits_{j=0}^{p-1}\sum\limits_{t=1}^e (ej+t)x_{ej+t}-2n)^2+\lambda\sum\limits_{t=1}^e y_{ej+t}^2 \\
&s.t.&   \bm g_j(\bm w_j)=((x_{ej+t}+y_{ej+t}-1)^2-y_{ej+t+n},\\
&& x_{ej+t}^2+(y_{ej+t}-1)^2-y_{ej+t+n},y_{ej+t}^2-y_{ej+t})=0,\ t=1,2,\cdots,e.
\end{eqnarray*}
For $j=0,1,2,\cdots,p-1$, the augmented Lagrange penalty optimizations of (Ex4.3)$_j$ are defined by
 \begin{eqnarray*}
\mbox{(Ex4.3)}_j &\min& \ \ A_j(\bm w_j,\bm \alpha_j,\sigma) = (\sum\limits_{j=0}^{p-1}\sum\limits_{t=1}^e (ej+t)x_{ej+t}-2n)^2+\lambda \sum\limits_{t=1}^e y_{ej+t}^2+\\ && \sum\limits_{t=1}^e [\alpha_{1t}((x_{ej+t}+y_{ej+t}-1)^2-y_{ej+t+n})+\alpha_{2t}(x_{ej+t}^2+(y_{ej+t}-1)^2\\ && -y_{ej+t+n}) +\alpha_{3}(y_{ej+t}^2-y_{ej+t})+\sigma(((x_{ej+t}+y_{ej+t}-1)^2-y_{ej+t+n})^2\\ && +(x_{ej+t}^2+(y_{ej+t}-1)^2-y_{ej+t+n})^2+\alpha_{3}(y_{ej+t}^2-y_{ej+t})^2)] \\
&s.t.&                (x_{ej+t},y_{ej+t},y_{ej+t+n})\in R^3, \ t=1,2,\cdots,e.
\end{eqnarray*}
By Algorithm 1,let the starting parameters $\epsilon=10^{-4},\sigma_{1}=5,N=10, \bm \alpha_j=(0, 0, \cdots,0)$ and $\bm w_j^0=(0,0,0,\cdots,0)$ be taken. When $e=5,p=1,2,6,10,20,$ $100,200$,$\lambda=1,10,100,500,1000$, value of
0-norm $\|\bm x^k\|_0$ is obtained by Algorithm 1 in Table. Numerical results show that an approximate sparse optimal solution is obtained in Table 1. When $\lambda$ is larger, value of
0-norm $\|\bm x^k\|_0$ is smaller at the approximate sparse optimal solution.

\begin{table}[ht]
%\small{
\begin{center}
\caption{ Value of $\|\bm x^k\|_0$  obtained by Algorithm 1 when $e=5$. }
\vspace{.001 in}
\begin{tabular}{|cccccc|}
 \hline
$n=5p$& $\lambda=1$ & $\lambda=10$ & $\lambda=100$& $\lambda=500$ &$\lambda=1000$\\
\hline
5  & 1 & 1  & 0 & 0 & 0\\
10 & 2 & 2 & 0 & 0& 0\\
30 &15  & 10 & 10  & 3 & 1 \\
50 & 22 &17 &18 &16 &10 \\
100 &54 &37 & 39 & 42 & 42 \\
500 &325 &192 & 208 & 227 & 208\\
1000 &611&382 & 427 & 471  & 421\\
\hline
\end{tabular}
\end{center}
\end{table}

{\footnotesize
\begin{table}[ht]
\begin{center}
\caption{ Numerical results  obtained by Algorithm 1 when $e=5$ and $n=100$. }
\vspace{.001 in}
\begin{tabular}{cccccccccccc}
 \hline
$\lambda$ &10 & 100& 1000&2000&3000&5000&8000&10000&15000&20000&21000\\
$\|\bm x^k\|_0$& 37&39&42&40&31&11&7&5&2&1&0\\
\hline
\end{tabular}
\end{center}
\end{table}
}
\end{example}

Numerical experiments show if the problem (CNP) is not decomposable, we cannot obtain an approximate sparse optimal solution by Algorithm 1 for $p=1$, $\lambda\in [1,10]$ and $n>15$. As the scale of the problem increases, the approximate spear solution may be obtained by using the decomposition algorithm (CNP). When $\lambda$ is greater, value of
0-norm $\|\bm x^k\|_0$ is smaller as shown in Table 2 when $n=e\times p=5\times 20=100$.

\begin{example} Consider nonconvex optimization problem  (Problem 5 in \cite{Bagirov},Problem 73 on Page 282 in \cite{Bagirov1} ):
 \begin{eqnarray*}
\mbox{(Ex4.4)} &\min& f_n(\bm x)=n\max\{|x_i|: i=1,2,\cdots,n\}-\sum\limits_{i=1}^n |x_i| \\
&s.t.& \bm x\in R^n.
\end{eqnarray*}
An optimal solution to (Ex4.4) is $\bm x^*=(\pm\alpha,\pm\alpha,\cdots,\pm\alpha)^\top$  with $f(\bm x^*)=0$ for $\alpha \in R^1$ in \cite{Bagirov1}. Let $\bm x\in R^n,\bm y\in R^{2n+1}$. A convertible nonconvex form of $f$ is defined by
 \begin{eqnarray*}
[n y_{2n+1}-\sum\limits_{i=1}^n y_i:
y_i^2-y_{i+n},x_i^2-y_{i+n},\ i=1,2,\cdots,n],
\end{eqnarray*}
where $\bm x\in S_1=R^n, \bm y\in S_2=\{\bm y\mid -y_i\leq 0,y_i-y_{2n+1}\leq 0,\ i=1,2,\cdots,n\}.$
Let $n=ep$, $\bm w_0=(x_1,x_2,\cdots,x_e,y_1,y_2,\cdots,y_e,y_{1+n},y_{2+n},\cdots,$
$y_{e+n},y_{2n+1})$, $\bm w_{j-1}=(x_{ej+1},x_{ej+2},\cdots,x_{ej+e},y_{ej+1},$ $y_{ej+2},\cdots, y_{ej+e},y_{ej+1+n},y_{ej+2+n},\cdots,y_{ej+e+n})$,
$j=1,2,\cdots,p-1$. So,for $j=0,1,2,\cdots,p-1$, $p$ subproblems are solved by:
 \begin{eqnarray*}
\mbox{(Ex4.4)}_j &\min& g(\bm w_0,\bm w_1,\cdots,\bm w_{p-1})=n y_{2n+1}-\sum\limits_{t=1}^e y_{et+j}\\
&s.t.&   \bm g_j(\bm w_j)=(y_{ej+t}^2-y_{ej+t+n},x_{ej+t}^2-y_{2n+1})=0,\ t=1,2,\cdots,e,\\
  &&             -y_{ej+t}\leq 0,y_{ej+t}^2-y_{ej+t+n}\leq 0,\ t=1,2,\cdots,e.
\end{eqnarray*}
For $j=0,1,2,\cdots,p-1$, the augmented Lagrange penalty optimizations of (Ex4.4)$_j$ are defined by
 \begin{eqnarray*}
\mbox{(Ex4.4)}_j &\min\ \ A_j(\bm w_j,\bm \alpha_j,\sigma) =& n y_{2n+1}-\sum\limits_{t=1}^e y_{et+j}+\sum\limits_{t=1}^e [\alpha_{1t}(y_{ej+t}^2-y_{ej+t+n})\\ && +\alpha_{2t}(x_{ej+t}^2-y_{2n+1})+\sigma((y_{ej+t}^2-y_{ej+t+n})^2\\
&& +(x_{ej+t}^2-y_{2n+1})^2+\max\{y_{ej+t}^2-y_{ej+t+n},0\}^2)] \\
&s.t.&                -y_{ej+t}\leq 0, \ t=1,2,\cdots,e.
\end{eqnarray*}
By Algorithm 1, let the starting parameters $\epsilon=10^{-4},\sigma_{1}=5,N=10,\bm \alpha_j=(0, 0, \cdots,0)$ and
$\bm w_j^0=(1,2,3,\cdots,3e+1)$ be taken. So, numerical results in Table 3 and Table 4 are obtained for $e=3,5$. Through Algorithm 1, an approximate solution can be obtained at iteration steps 2-4 . The numerical results show that the running time is a multiple of the decomposition scale $p$. Through numerical experiments, it is appropriate to take $e$ from 2 to 6 in this example. When $e$ exceeds 6, the larger $n$, the worse the approximate solution. However, the smaller the $e$, the longer the execution time. The scale of (Ex4.4)  here is much larger than those of (Ex4.4) in literature \cite{Bagirov}. If fact, the example shows that the the solution $\bm x^*$ to (Ex4.4) is composed of $\bm x^k$ repeated $e$ times. For example, because $\bm x^k=(3.1448,3.1448,3.1448, 3.1448,3.1448)$ is an optimal solution to (Ex4.4) at $n=5$ by Algorithm 1. Then, $(3.1448,3.1448,3.1448, 3.1448,$ $ 3.1448)$ extends to $\bm x^k=(3.1448,3.1448,3.1448, 3.1448,3.1448,3.1448,3.1448,$ $3.1448, 3.1448,3.1448)$ which is an optimal solution to (Ex4.4) at $n=10$. When $n$ is very large, it is easy to obtain its optimal solution.

\begin{table}[ht]
%\small{
\begin{center}
\caption{ Numerical results  obtained by Algorithm 1 when $e=5$. }
\vspace{.001 in}
\begin{tabular}{cccc}
 \hline
$n$& $k$  & $(x_1^k,x_2^k,\cdots,x_n^k)$ & Running time \\
\hline
5& 2&(3.1448,3.1448,3.1448, 3.1448,3.1448)&2.9155s\\
10& 2&$(3.1448,3.1448,3.1448,\cdots,3.1448)$&4.3817s\\
50& 2&$(3.1448,3.1448,3.1448,\cdots,3.1448)$&18.1406s\\
250& 3&$(3.3481,3.3481,3.3481,\cdots,3.3481)$&152.7852s\\
500& 4&$(3.3484,3.3484,3.3484,\cdots,-3.3484)$&401.494226s\\
1000& 4&$(3.3484,3.3484,3.3484,\cdots,3.3484)$&855.879704s\\
\hline
\end{tabular}
\end{center}
\end{table}
\begin{table}[ht]
%\small{
\begin{center}
\caption{ Numerical results  obtained by Algorithm 1 when $e=3$. }
\vspace{.001 in}
\begin{tabular}{cccc}
 \hline
$n$& $k$  & $(x_1^k,x_2^k,\cdots,x_n^k)$ & Running time \\
\hline
6& 2&$(3.4322,3.4322,3.4322,\cdots,3.4322)$&3.1981s\\
30& 2&$(3.4322,3.4322,3.4322,\cdots,3.4322)$&12.5217s\\
90& 2&$(3.4322,3.4322,3.4322,\cdots,3.4322)$&45.0358s\\
300& 3&$(3.4262,3.4262,3.4262,\cdots,3.4262)$& 215.3782s\\
600& 3&$(3.4262,3.4262,3.4262,\cdots,3.4262)$& 457.0393s\\
1500& 3&$(3.4262,3.4262,3.4262,\cdots,3.4262)$&1102.5697s\\
\hline
\end{tabular}
\end{center}
\end{table}
\end{example}

\begin{example} Consider a nonconvex optimization problem  (Problem 64 on Page 280 in \cite{Bagirov1} ):
 \begin{eqnarray*}
\mbox{(Ex4.5)} &\min& f_n(\bm x)=\sum\limits_{i=1}^{n-1}(- x_i+2(x_i^2+x_{i+1}^2-1)+1.75|x_i^2+x_{i+1}^2-1|) \\
&s.t.& \bm x\in R^n.
\end{eqnarray*}
 Let $\bm x\in R^n,\bm y\in R^{2n+1}$. A convertible nonconvex form of $f$ is defined by
 \begin{eqnarray*}
[\sum\limits_{i=1}^{n-1}(-x_i+ 2y_i+1.75y_{i+n})&:&
x_i^2+x_{i+1}^2-1-y_i=0,y_{i+n}^2-y_{i+2n}=0,\\ && y_i^2-y_{i+2n}=0\ i=1,2,\cdots,n-1],
\end{eqnarray*}
where $\bm x\in S_1=R^n, \bm y\in S_2=\{\bm y\mid y_i\geq -1,y_{i+n},y_{i+2n}\geq 0,\ i=1,2,\cdots,n\}.$
Let $n=ep$, $\bm w_{j}=(x_{ej+1},x_{ej+2},\cdots,x_{ej+e},y_{ej+1},y_{ej+2},\cdots, y_{ej+e},y_{ej+1+n},$ $ y_{ej+2+n},\cdots,y_{ej+e+n},y_{ej+1+2n},y_{ej+2+2n},$ $\cdots,y_{ej+e+2n})$,
$j=0,1,2,\cdots,p-1$. So, for $j=0,1,2,\cdots,p-1$, $p$ subproblems are solved by:
 \begin{eqnarray*}
\mbox{(Ex4.5)}_j &\min& g(\bm w_0,\bm w_1,\cdots,\bm w_{p-1})=\sum\limits_{t=1}^e (-x_{ej+t}+ 2y_{ej+t}+1.75y_{ej+t+n}) \\
&s.t.&   \bm g_j(\bm w_j)=(x_{ej+t}^2+x_{ej+t+1}^2-1-y_{ej+t},y_{ej+t+n}^2-y_{ej+t+2n},\\
&& y_{ej+t}^2-y_{ej+t+2n})=0,\ t=1,2,\cdots,e,\\
  &&             -y_{ej+t}\leq 1,-y_{ej+t+n},-y_{ej+t+2n}\leq 0,\ t=1,2,\cdots,e.
\end{eqnarray*}
Since $\{\bm w_{j-1}\}\cap \{\bm w_{j}\}=\{x_{ej+1}\}$, $f_n(\bm x)$ is not decomposable. If $\bm w_{j-1}^*(j>0)$ is an optimal solution to (EX4.5)$_{j-1}$, a term "$\sigma(x_{ej+1}^*-x_{ej+1})^2$" is added to (EX4.5)$_j$, where for $j=0,1,2,\cdots,p-1$, the augmented Lagrange penalty optimizations of (EX4.5)$_j$ are defined by
 \begin{eqnarray*}
\mbox{(Ex4.5)}_j &\min\ \ A_j(\bm w_j,\bm \alpha_j,\sigma) =& \sum\limits_{t=1}^e (-x_{ej+t}+ 2y_{ej+t}+1.75y_{ej+t+n})+\\
&& \sum\limits_{t=1}^e [\alpha_{1t}(x_{ej+t}^2+x_{ej+t+1}^2-1-y_{ej+t})+\\
&& \alpha_{2t}(y_{ej+t+n}^2-y_{ej+t+2n})+\alpha_{3t}(y_{ej+t}^2-y_{ej+t+2n}))\\
&& +\sigma((x_{ej+t}^2+x_{ej+t+1}^2-1-y_{ej+t})^2+\\
&& (y_{ej+t+n}^2-y_{ej+t+2n})^2+(y_{ej+t}^2-y_{ej+t+2n})^2)] \\
&& +\sigma(x_{ej+1}^*-x_{ej+1})^2\\
&s.t.&               -y_{ej+t}\leq 1,-y_{ej+t+n},-y_{ej+t+2n}\leq 0, \ t=1,2,\cdots,e.
\end{eqnarray*}
Especially, when $j=0$, "$\sigma(x_{1}^*-x_{1})^2$" is deleted from the above (EX4.5)$_j$.
When $n=50, 200$ and $1000$, their approximate objective values $f(\bm x^k)$ are $-34.795,-140.86$ and $-706.55$ respectively at starting point $x^1=(1,1,\cdots,1)^\top$ in literature \cite{Bagirov1}.
By Algorithm 1,let the starting parameters $\epsilon=10^{-4},\sigma_{1}=5,N=100,\bm \alpha_j=(0, 0, \cdots,0)^\top$ and
$\bm w_j^0=(1,1,1,\cdots,1)^\top$ be taken. So, at iteration step 1, numerical results are obtained by Algorithm 1 as shown in Table 5, 6 and 7 respectively when $n=50,200$ and $1000$. In Table 5,6 and 7, the best objective value of (Ex4.5) are obtained when $n=2p$. It is worth noting that the solutions are effective in Table 5-7, because (Ex4.5) is an unconstrained optimization problem.

\begin{table}[ht]
%\small{
\begin{center}
\caption{ Numerical results  obtained by Algorithm 1 when $n=50$. }
\vspace{.001 in}
\begin{tabular}{ccccc}
 \hline
$e$& $p$  & $f_n(\bm x^k) $ & $\|\bm g(\bm x^k,\bm y^k)\|$ &Running time(s)\\
\hline
50 & 1 & -38.2547  & 0.0593 & 4.2525 \\
25 & 2 & -36.8531  & 0.0472 & 5.6183 \\
10 & 5 & -39.4672  & 0.0001 &10.8266 \\
5 & 10 & -35.8447  & 0.0000 &24.6390\\
2 & 25 & -40.8154  & 0.0000 &47.2351\\
\hline
\end{tabular}
\end{center}
\end{table}

\begin{table}[ht]
%\small{
\begin{center}
\caption{ Numerical results  obtained by Algorithm 1 when $n=200$. }
\vspace{.001 in}
\begin{tabular}{ccccc}
 \hline
$e$& $p$  & $f_n(\bm x^k) $ & $\|\bm g(\bm x^k,\bm y^k)\|$ &Running time(s)\\
\hline
%200 & 1 & -214.1186  & 0.77871 & 7.2784\\
%100 & 2 & -210.8638  & 0.36372 &10.4066 \\
20 & 10 & -159.4185  & 0.04810 &25.5085\\
10 & 20 & -159.4124  & 0.00101 &37.1376\\
5 & 40 & -145.9025  & 0.00002 &82.4633\\
4 & 50 & -156.7485  & 0.00000 &93.6619\\
2 & 100 & -166.0306  & 0.00000 &151.7478\\
\hline
\end{tabular}
\end{center}
\end{table}

\begin{table}[ht]
%\small{
\begin{center}
\caption{ Numerical results  obtained by Algorithm 1 when $n=1000$. }
\vspace{.001 in}
\begin{tabular}{ccccc}
 \hline
$e$& $p$  & $f_n(\bm x^k) $ & $\|\bm g(\bm x^k,\bm y^k)\|$ &Running time(s)\\
\hline
10 & 100 & -805.0714  & 0.00108 &271.2822\\
%10 & 100 & -801.5056  & 0.00111 &227.9784
5 & 200 & -726.9435  & 0.00002 &566.6910\\
%5 & 200 & -727.4786  & 0.00003 &476.3201\\
4 & 250 & -783.3616  & 0.00001 &589.5969\\
2 & 500 & -833.6333  & 0.00002 &717.9175\\
\hline
\end{tabular}
\end{center}
\end{table}

By Algorithm 1, a solution $\tilde{\bm x}=(0.6897,0.4153,0.4965, 0.5386,0.4707)$ with its objective value $f_5(\tilde{\bm x})=-3.0766$ is obtained at $n=5$ and $p=1$. And another solution $\hat{\bm x}=(0.6772,0.49930)$ with its objective value $f_2(\hat{\bm x})=-0.8749$ is obtained at $n=2$ and $p=1$. In Table 8,
the solution $\bm x^k$ of the first three columns is composed of $\tilde{\bm x}$ repeated $p$ times, such as
$\bm x^k=(\tilde{\bm x},\tilde{\bm x},\tilde{\bm x},\tilde{\bm x},\tilde{\bm x},\tilde{\bm x},\tilde{\bm x},\tilde{\bm x},\tilde{\bm x},\tilde{\bm x})$ at $e=5$ and $p=10$ in line 1, where $f_n(\bm x^k)$ is the objective value of $\bm x^k$.
In Table 8, the solution $\bm x^k$ of the from 4 to 6 columns is composed of $\hat{\bm x}$ repeated $p$ times.
In Table 8, the values of the last six columns are obtained by Algorithm 1 (See Table 5,6,7).
Numerical results show that for large-scale unconstrained optimization problems, a better solution can be obtained directly by Algorithm 1 by solving  small-scale subproblems of it  when the structure of all the subproblems are similar, i.e. constraint structure of all equality subproblems is the same as in Example 4.3,4.4 and 4.5. This decomposable method is effective in examples 4.4 and 4.5.

\begin{table}[ht]
{\tiny
\begin{center}
\caption{ Numerical results  obtained by Algorithm 1. }
\vspace{.001 in}
\begin{tabular}{|ccc|ccc||ccc|ccc|}
 \hline
$e$& $p$  & $f_n(\bm x^k) $ & $e$& $p$  & $f_n(\bm x^k)$&$e$& $p$  & $f_n(\bm x^k) $ & $e$& $p$  & $f_n(\bm x^k)$\\
\hline
5 & 10 & -37.8238  & 2 & 25 & -41.1109& 5 & 10 & -35.8447 & 2 & 25 & -40.8154  \\
5 & 40 & -153.6478  & 2 & 100 & -166.8484 &5 & 40 & -145.9025  & 2 & 100 & -166.0306\\
5 & 200 & -771.3758  & 2 & 500 & -837.4484& 5 & 200 & -726.9435 & 2 & 500 & -833.6333\\
5 & 1000 & -3860.0158& 2 & 2500 & -4190.4484& 5 & 1000 & -3647.9851 & 2 & 2500 & -4171.0713\\
5 & 2000 & -7720.8158& 2 & 5000 & -8381.6984& 5 & 2000 & -7295.0801 &2 & 5000 & -8343.5399  \\
\hline
\end{tabular}
\end{center}
}
\end{table}

%5 & 2000 & -7295.0801  & 0.00025 &3927.0610 "
%2 & 5000 & -8343.5399  & 0.00017 &6966.8208 "\"

\end{example}

Algorithm 1 shows that when these parameters $\epsilon>0,\sigma_{1}>0,N>1, (\bm x_j^0,\bm y_j^0)$ $(j=1,2,\cdots,p), \bm \alpha_j^1$$(j=1,2,\cdots,p)$ are properly selected, it may obtain an approximate optimal solution to an unconstrained optimization problem with CN function. And above examples show that the scale problem of (CNP) can be avoided and reduced by solving the subproblems in Algorithm 1.\\

{\bf Note:} If (CNP) is a small-scale problem for $p=1$ or the number of variables in (CNP) is very small, Algorithm 1 may obtain an approximate solution to (CNP).

\section{Conclusion}

This paper solves three difficulties relating to unconstrained nonconvex or nonsmooth optimization problems. (1) An unconstrained, nonconvex and nonsmooth optimization problem is transformed into a constrained optimization problem, where its objective function and constrained functions are convex and smooth.  A new concept - weak uniform or decomposable CN optimization - is proposed, which covers many nonconvex nonsmooth functions, even discontinuous nonconvex functions. (2) The optimality conditions of the global optimal solution to unconstrained nonconvex or nonsmooth problems are obtained. The sufficient conditions of the global optimal solution are proved for the weak uniform or decomposable CN optimization problems. (3) A decomposable algorithm for unconstrained nonconvex or nonsmooth optimization problem is proposed, based on the augmented Lagrange penalty function of (CNP).
The results of  numerical examples show that Algorithm 1 may obtain an approximate global optimal solution after properly selecting the initial parameters. Another advantage of the algorithm is that it does not need to use subgradient and smoothing techniques such that Matlab-the optimization software is directly usable, which makes this method easy for engineers to use.

This paper provides a new method for solving nonconvex unconstrained optimization problems, which shows its potential importance  in many application fields. There are at least three directions worthy of further study in terms of weak uniform CN optimization problems (CNP):

(1) decomposable Newton algorithm or decomposable SQP algorithm for (CNP),

(2) Lagrangian multiplier alternating algorithm for (CNP),

(3) some special structures of (CNP), for example when $f$ is a weak uniform CN function with $f=[g,g_1,g_2,\cdots,g_r]$, where  $g,g_1,g_2,\cdots,g_r$ are quadratic functions.

\baselineskip 18pt

{\bf Acknowledgements}

This work is supported by the National Natural Science Foundation of China(No.11871434).

%\end{acknowledgements}

{\bf Data Availability Statement}

The author confirms that all data generated or analysed during this study are included in the published article \cite{Bagirov,Bagirov1,Chen1,Chen2}.

%The authors  would like to express their gratitudes to anonymous referees' detailed comments and remarks that help us improve our presentation of this paper considerably.
\bibliographystyle{siamplain}
%\begin{thebibliography}{plain}

\end{document}